\documentclass[english,11pt,leqno,twoside]{amsart}
\usepackage[T1]{fontenc}
\usepackage{amsmath}
\usepackage{amsthm}
\usepackage{amssymb}
\usepackage{amsfonts}
\usepackage[left=2cm,right=2cm]{geometry}
\usepackage{enumitem}

\usepackage[noadjust]{cite}
\usepackage{lmodern}
\usepackage{hyperref}
\usepackage[hyperpageref]{backref}
\usepackage{fixmath}
\usepackage{todonotes}
\usepackage{csquotes}

\usepackage{fancyhdr}

\numberwithin{equation}{section}
\newtheorem{Th}{Theorem}[section]
\newtheorem{Prop}[Th]{Proposition}
\newtheorem{Lem}[Th]{Lemma}
\newtheorem{Cor}[Th]{Corollary}
\newtheorem{Def}[Th]{Definition}
\newtheorem{Rem}[Th]{Remark}

\newcommand{\N}{\mathbb{N}}
\newcommand{\R}{\mathbb{R}}

\newcommand{\cI}{\mathcal{I}}
\newcommand{\cJ}{\mathcal{J}}

\newcommand{\cS}{{\mathcal S}}

\newcommand{\loc}{\mathrm{loc}}

\newcommand{\dist}{\mathrm{dist}}
\newcommand{\les}{\lesssim}

\newcommand{\eps}{\varepsilon}
\renewcommand{\phi}{\varphi}
 \def\dd{\, {\rm d}}
\newcommand{\tf}{\widetilde f}
\newcommand{\tu}{\widetilde u}
\newcommand{\tv}{\widetilde v}
\newcommand{\tw}{\widetilde w}
\newcommand{\tX}{\widetilde X}
\newcommand{\tGamma}{\widetilde\Gamma}

\newcommand{\intd}{\int_{\R^2}}
\def\e{{\rm e}}
\newcommand{\ff}{\mathfrak f}

\title{Nonlinear Schr\"odinger-Poisson systems in dimension two:\\ the zero mass case}
\author{Federico Bernini, Giulio Romani, Cristina Tarsi}

 \address[Federico Bernini]{\newline\indent Dipartimento di Matematica e Applicazioni,
         \newline\indent Università degli Studi di Milano-Bicocca,
         \newline\indent Via Cozzi 55 - 20125 Milano (Italy)}
         \email{\href{mailto:federico.bernini@unimib.it}{federico.bernini@unimib.it}}

 \address[Giulio Romani]{
 	\newline\indent Dipartimento di Scienze Matematiche, Informatiche e Fisiche,
 	\newline\indent Universit\`{a} degli Studi di Udine
 	\newline\indent Via delle Scienze 206 - 33100 Udine, Italy}
     \email{\href{mailto:giulio.romani@uniud.it}{giulio.romani@uniud.it}}

 \address[Cristina Tarsi]{\newline\indent Dipartimento di Matematica,
	\newline\indent Università degli Studi di Milano,
	\newline\indent Via C. Saldini, 50 - 20133 Milan, Italy}
\email{\href{mailto:cristina.tarsi@unimi.it}{cristina.tarsi@unimi.it}}

\allowdisplaybreaks

\begin{document}

\begin{abstract}
We provide an existence result for a Schr\"odinger-Poisson system in gradient form, set in the whole plane, in the case of zero mass. Since the setting is limiting for the Sobolev embedding, we admit nonlinearities with subcritical or critical growth in the sense of Trudinger-Moser. In particular, the absence of the mass term requires a nonstandard functional framework, based on homogeneous Sobolev spaces. These features, combined with the logarithmic behaviour of the kernel of the Poisson equation, make the analysis delicate, since standard variational tools cannot be applied. The system is solved by considering the corresponding logarithmic Choquard equation. The existence of a mountain pass-type solution is established by means of a careful analysis of appropriate Cerami sequences, whose boundedness is ensured through a nonstandard variational method, suggested by the subtle nature of the functional geometry involved. As a key tool in our estimates, we also introduce a logarithmic weighted Trudinger–Moser inequality, along with a related Cao-type inequality, both of which hold in our functional setting and are, we believe, of independent interest.
\end{abstract}


\maketitle
\section{Introduction}
We aim at investigating existence of positive solutions of the planar Schr\"odinger-Poisson system in gradient form given by 
\begin{equation}\label{SP:system}
	\begin{cases}
		-\Delta u = \Phi f(u) &\text{ in } \R^2,\\
		-\Delta\Phi = 2\pi F(u) &\text{ in } \R^2,
	\end{cases}
\end{equation}
where $f$ is a positive continuous nonlinearity with subcritical or critical growth in the sense of Trudinger-Moser, and $F(t):=\int_0^tf(s)\dd s$. The main goal is to face the combined difficulties of working in the limiting setting of the Sobolev embeddings and the fact that in the first equation of \eqref{SP:system} the mass term is missing. This makes the problem challenging, not only for the variational approach, but also in the choice of a nonstandard functional framework.
\vskip0.2truecm
Schr\"odinger-Poisson systems of the form
\begin{equation}\label{SP:system-RN-V}
	\begin{cases}
		-\Delta u +V(x)u= \Phi f(u) &\text{ in } \R^N,\\
		-\Delta\Phi = F(u) &\text{ in } \R^N,
	\end{cases}
\end{equation}
with a potential $V$ which is usually positive, are of great importance in several fields of physics, since they serve as models for the interaction of two identically charged particles in electromagnetism, as well as for the self-interaction of the wave function with its own gravitational field in quantum mechanics; they appear also in the Hartree theory for crystals, and in astrophysics in the study of selfgravitating boson stars; for the physics background we refer to \cite{BF,LRZ} and to the references therein. From a mathematical point of view, they are interesting since they can be analysed by variational techniques. Indeed, one may (formally) solve the Poisson equation in \eqref{SP:system} by means of the Riesz kernel
$$K_N(x):=\begin{cases}
	\frac{C_N}{|x|^{N-2}}&\text{ if } N\geq3\,,\\
	\frac1{2\pi}\ln\frac1{|x|} &\text{ if } N=2\,,
\end{cases}$$
where $C_N$ is an explicit positive constant, and consider
\[
\Phi_u(x):=\left(K_N\ast F(u)\right)(x)=\int_{\R^N}K_N(x-y)F(u(y))\dd y\,.
\]
Substituting it into the first equation of \eqref{SP:system-RN-V}, one may rewrite the system as a Choquard equation, that is an integro-differential equation of Schr\"odinger type with a convolutive right-hand side:
\begin{equation}\label{SP:equation-RN-V}
	-\Delta u(x) + V(x)u = \left(K_N\ast F(u)\right)f(u)\quad\text{in}\ \R^N.
\end{equation}
Besides the evident advantage of the reduction to a single equation, since the system \eqref{SP:system-RN-V} is of gradient type, if $N\geq3$ the equation \eqref{SP:equation-RN-V} is variational in the Sobolev space $H^1(\R^N)$ thanks to the Hardy-Littlewood-Sobolev inequality (Proposition \ref{HLS} below), provided suitable polynomial growth conditions on $f$ are fulfilled. In this respect, there are a huge number of works about Choquard-type equations, especially from the last decades, and we refer to the seminal works \cite{MV,MV2,CZ,CVZ} and the references therein. The planar case $N=2$ is more delicate because of the interplay between the logarithmic behaviour of the Riesz kernel, and the exponential maximal growth of the nonlinearities, due to the Poho\v zaev-Trudinger-Moser inequality in the full space proved by Ruf \cite{Ruf}, see also \cite{Cao}. The first attempt of considering this case is to be referred to Stubbe \cite{Stubbe2008}, later on formalised by Cingolani and Weth \cite{CingolaniWeth2017,CW2}: they set the problem in a constrained space which takes into account in the seminorm a contribution of the logarithmic kernel. This analysis, which is peculiar for the case of a linear coupling in the system, namely $f(u)=u$ in \eqref{SP:system-RN-V}, was then extended for the general case of a nonlinearity with critical exponential growth in \cite{CassaniTarsi2021}. Taking indeed into account the behaviour at $0$ of the nonlinearity, and by means of a log-weighted version of the Poho\v zaev-Trudinger inequality, a proper functional setting was found, in which the functional associated to \eqref{SP:equation-RN-V} turns out to be well-defined. We also point out that the sharp version of such inequality has been recently obtained in \cite{Tarsi}. The approach in \cite{CassaniTarsi2021} was then generalised for Choquard equations with weights in \cite{CMR} and for quasilinear Schr\"odinger-Poisson systems in \cite{BCT}. A different approach was recently proposed in \cite{LRTZ}: here, instead, the underlying functional space remains $H^1(\R^2)$, while the logarithmic kernel is uniformly approximated by polynomial kernels. For further developments of this method we refer to \cite{CDL,CLR,CLR2}, also in quasilinear fractional contexts. 
\vskip0.2truecm
All the above works deal with Choquard equations of the form \eqref{SP:equation-RN-V}, where $V$ is a nontrivial potential. The special case of an identically zero potential, the so-called ``zero-mass case'', emerges in some physical context, e.g. in the nonabelian gauge theory of particle physics, such as the study of the Yang-Mills equation, see \cite{Gidas}. The main difference with the ``mass-case'' lies in the natural framework in which such problems are studied, namely the homogeneous Sobolev spaces $D^{1,2}_0(\R^N)$ defined as the completion of $C^\infty_c(\R^N)$ with respect of the $L^2$-norm of the gradient. In the higher-dimensional case $N\geq3$, in light of the critical Sobolev embedding into $L^{2^*}(\R^N)$ with $2^*=\frac{2N}{N-2}$, such setting is appropriate for both Schr\"odinger and Choquard equations of the kind
\begin{equation}\label{SP:equation-RN-0}
	-\Delta u = \left(K_N\ast F(u)\right)f(u)\quad\text{in}\ \R^N,
\end{equation}
as shown e.g. in \cite{BL,AP,ASM,AY}. However, if $N=2$, not only any embedding into Lebesgue spaces is out of reach, but also one cannot distinguish in $D^{1,2}_0(\R^2)$ between functions which just differ by constants. In this framework only few results are available for Schr\"odinger equations and Choquard equations of the kind \eqref{SP:equation-RN-0} with $K_2=\frac1{2\pi}\ln\frac1{|\cdot|}$. On the one hand, one may try to adjust the operator so that the corresponding natural space recovers good embedding properties in Lebesgue and Orlicz spaces: this is the strategy used in \cite{dAC,CFFM} for the Schr\"odinger case, and in \cite{Romani1,Romani2,dACS} for the Choquard case with both polynomial and logarithmic kernels. On the other hand, for equations with zero mass driven by the pure Laplacian, the results available in the literature \cite{ChenRadulescuWen,ChenTang,ChenShuTangWen2022} cover just the linear case $f(u)=u$ in \eqref{SP:equation-RN-0} (up to adding local nonlinearities) since the approach of \cite{Stubbe2008,CingolaniWeth2017} is pursued. However, as remarked in \cite{ChenShuTangWen2022},
\begin{displayquote}
	\textit{``it is really interesting to observe how remarkable the impact of the logarithmic integral kernel $\ln|x|$ is, because it allows us to establish much richer and better existence results than those available for other elliptic equations, in spite of its sign-changing and unbounded properties''}.
\end{displayquote}
Indeed, differently from the Schr\"odinger equations with zero mass, the presence of the logarithmic kernel combined with $f(u)=u$ allows to recover $H^1(\R^2)$ as suitable functional framework, thanks to a careful splitting in positive and negative part of the logarithm (see \eqref{logarithmic:equality} below).
\vskip0.2truecm
Main goal of this paper is to extend the existence results in the zero-mass case of \cite{ChenRadulescuWen,ChenTang,ChenShuTangWen2022} in the direction of \cite{CassaniTarsi2021}, that is to consider the general case of a nonlinear function $f$, and cover both the cases of subcritical and critical growth in the sense of the Trudinger-Moser inequality. To this aim, several difficulties need to be faced: first, the unusual functional setting, which does not appear in the above cited works because of the linear behaviour of $f$; then the analysis on Cerami sequences, which arise from the mountain pass geometry of the functional associated to \eqref{SP:equation}, is largely affected by the possibly exponential growth of the nonlinear terms and, in particular in the critical case, it is very delicate; eventually the final proof ``à la Lions'' of the existence theorem is also quite non-standard and rather technical. Finally, we derive from our results for the Choquard equation the corresponding for the Schr\"odinger-Poisson system \eqref{SP:system} in a suitable functional setting in the spirit of  \cite{BCT,Romani1}. We stress that this step is often neglected in the literature by just considering it as ``natural'': here it finds a rigorous justification.
\vskip0.2truecm
Before stating our main results, let us specify the growth conditions we are assuming on the nonlinearity. Throughout the paper, we suppose that $f\in C^1(\R)$, $f(s)>0$ as $s>0$, while $f(s)=0$ for $s\leq 0$; moreover, it satisfies:
\begin{itemize}
	\item[$(f_1)$] $f(s)\asymp s^{p-1}$ as $s\to 0$ for some $p>2\,$;
\end{itemize}
and either is subcritical or critical in the following sense:
\begin{itemize}
	\item [$(f_{sc})$] for any $\alpha>0$, $\displaystyle\lim_{s\to +\infty} f(s)/\e^{\alpha s^2}=0$ and for some $C>0$, $f(s)\geq Cs^{p-1}$ as $s\to+\infty\,$, with $p>2$;
	\item [$(f_c)$] there exists $\alpha_0>0$ such that
	$$\lim_{s\to+\infty}\frac{f(s)}{\e^{\alpha s^2}}=\begin{cases}
		0 &\text{ if } \alpha>\alpha_0\,,\\
		+\infty &\text{ in } \alpha<\alpha_0\,.
	\end{cases}$$
\end{itemize}
Under such conditions, and in particular in light of the behaviour near $0$ given by $(f_1)$, we are going to see that the functional setting in which it is convenient to look for solutions of the Choquard equation associated to the system \eqref{SP:system}, namely 
\begin{equation}\label{SP:equation}
	-\Delta u + \left(\ln|\cdot|\ast F(u)\right)f(u) = 0\quad\text{in}\ \R^2,
\end{equation}
is
$$D^{1,2}L^p_\omega(\R^2):=D^{1,2}_0(\R^2)\cap L^p(\R^2,\omega\dd x)\,,$$
with weight function $\omega(x):=\ln(b+|x|)$ and $b>1$. Note that this space, which corresponds to the intersection space $H^1L^p_\omega(\R^2)$ detected in \cite{CassaniTarsi2021} for the study of the same Choquard equation with positive mass, gathers all important features of our problem: the absence of mass since we are dealing with the homogeneous Sobolev space $D^{1,2}(\R^2)$, the nonlinear behaviour of $f$, and the logarithmic kernel in the weight. We will see that this space enjoys good embedding properties, in particular exponential nonlinearities are allowed. We are therefore led to the following definition:
\begin{Def}[Solution of \eqref{SP:equation}]\label{Def_Choquard}
	We say that $u\in D^{1,2}L^p_\omega(\R^2)$ is a \emph{weak solution} of \eqref{SP:equation} if
	\begin{equation*}
		\intd\nabla u \cdot \nabla\varphi\dd x=\intd\!\left(\intd\ln\frac1{|x-y|}F(u(y))\dd y\!\right)f(u(x))\varphi(x)\dd x
	\end{equation*}
	for all $\varphi\in D^{1,2}L^p_\omega(\R^2)$.
\end{Def}
Of course, in order to prove existence for \eqref{SP:equation}, we need some further assumptions on $f$, which are gathered here:
\begin{itemize}
	\item[$(f_2)$] there exist $C>\tau>0$ such that $\tau\leq\frac{F(s)f'(s)}{f(s)^2}\leq C$ for all $s>0\,$;
	\vskip0.1truecm
	\item[$(f_3)$] $\displaystyle\lim_{s\to +\infty}\frac{F(s)f'(s)}{f(s)^2}=1$, or equivalently $\displaystyle\lim_{s\to +\infty}\frac{\dd}{\dd s}\frac{F(s)}{f(s)}=0\,$;
	\vskip0.1truecm	
	\item[$(f_4)$] $\displaystyle\lim_{s\to +\infty}\frac{s^3f(s)F(s)}{\e^{2\alpha_0s^2}}\geq\beta>\mathcal V$, where $\mathcal V$ will be explicitly given in \eqref{V};
	\vskip0.1truecm
	\item[$(f_5)$] $f'(s)\asymp s^{p-2}$ as $s\to 0$ with $p>2$, and $f'$ satisfies ($f_c$) with the same $\alpha_0$ as $f$.	
\end{itemize}
\vskip0.1truecm
We postpone to Section \ref{Assumption:consequences} a detailed list of consequences of our assumptions. Here we just emphasise that $(f_2)$ implies an Ambrosetti-Rabinowitz condition and the monotonicity of $f$; $(f_3)$-$(f_4)$ will be used in the analysis of the boundedness of Cerami sequences when dealing with critical nonlinearities, the latter being related to the de Figueiredo-Miyagaki-Ruf condition in \cite{dFMR} and used to prove a fine upper bound for the mountain pass level in Section \ref{Section_Cerami}; $(f_5)$ is a mild condition about the growth at $\infty$ of the nonlinearity, which well agrees with both $(f_{sc})$-$(f_c)$, and will be exploited in the conclusive compactness argument in Section \ref{Section_existence}.
\vskip0.2truecm
Our main result reads as follows.
\begin{Th}[Existence for \eqref{SP:equation}]\label{existence:nontrivial:solution}
	Suppose ($f_1$), ($f_2$), ($f_5$) hold, and either
	\begin{itemize}
		\item[i)] $f$ is subcritical as in ($f_{sc}$),
	\end{itemize}
	or
	\begin{itemize}
		\item[ii)] $f$ is critical as in ($f_c$) and ($f_3$), ($f_4$) are fulfilled.
	\end{itemize}
	Then there exists a positive solution to \eqref{SP:equation} in the sense of Definition \ref{Def_Choquard}.
\end{Th}
Once we have found a weak solution of the logarithmic Choquard equation \eqref{SP:equation}, we can go back to the original Schr\"odinger-Poisson system. First, we need a precise meaning of solution for \eqref{SP:system}.
The weighted Lebesgue space $L_s(\R^2)$, $s>0$, is defined as
$$
L_s(\R^2):=\Big\{u\in L^1_{loc}(\R^2)\,\Big|\,\int_{\R^2}\frac{|u(x)|}{1+|x|^{2+2s}}\dd x<+\infty\Big\}\,.
$$
\begin{Def}\label{sol_Poisson}
	For $\ff\in\mathcal S'(\R^2)$ we say that a function $\phi\in L_1(\R^2)$ is a \emph{solution} of the linear Poisson equation $-\Delta\Phi=\ff$ in $\R^2$ if
	$$
	\intd\Phi\,(-\Delta\varphi)=\langle\ff,\varphi\rangle\qquad\mbox{for all}\ \,\varphi\in\mathcal S(\R^2)\,.
	$$
\end{Def}
\begin{Def}[Solution of \eqref{SP:system}]\label{sol_SP}
	We say that $(u,\Phi)$ is a \emph{weak solution} of the Schr\"odinger-Poisson system \eqref{SP:system} if
	\begin{equation*}
		\intd\nabla u\cdot\nabla\varphi\dd x=\intd\Phi f(u)\varphi\dd x
	\end{equation*}
	for all $\varphi\in D^{1,2}L^p_\omega(\R^2)$, and $\Phi$ solves $-\Delta\Phi=2\pi F(u)$ in $\R^2$ in the sense of Definition \ref{sol_Poisson}.
\end{Def}
\begin{Th}[Existence for \eqref{SP:system}]\label{Thm_SP}
	Under the conditions of Theorem \ref{existence:nontrivial:solution}, the Schr\"odinger-Poisson system \eqref{SP:system} possesses a solution $(u,\Phi)\in D^{1,2}L^p_\omega(\R^2)\times L_s(\R^2)$ for all $s>0$, such that $u$ is positive and $\Phi=\Phi_u:=\ln\frac1{|\cdot|}\ast F(u)$.
\end{Th}
\begin{Rem}\label{Rem_pos_sol}
	It is worth to point out that:
	\begin{enumerate}
		\item this work can be seen as an extension to the zero mass-case of the results in \cite{CassaniTarsi2021}, to the general case of a nonlinearity with possibly exponential growth of the results in \cite{ChenRadulescuWen,ChenTang,ChenShuTangWen2022}, and to the pure Laplacian case to those in \cite{Romani1};
		\item it is sufficient to prove the existence of a nonnegative nontrivial solution of \eqref{SP:equation} in order to retrieve its positivity by the strong maximum principle for semilinear equations, see e.g. \cite[Theorem 11.1]{PucciSerrin}.
	\end{enumerate}
\end{Rem}
\paragraph{\textbf{Overview}} In Section \ref{Section_spaces} we describe the functional framework in which it is convenient to set our problem, and prove Trudinger-Moser- and Cao-type inequalities in weighted spaces, which will be fundamental tools in our arguments; moreover, we discuss our assumptions and collect some useful results. The variational framework is then described in Section \ref{Section_Var}, where we show the mountain pass geometry for the energy functional, while the existence of special Cerami sequences, and their boundedness is detailed in Section \ref{Section_Cerami}; we stress that these arguments turn out to be a delicate matter. After some careful mountain pass estimates, the proof of the existence for the log-Choquard equation \eqref{SP:equation} is given in Section \ref{Section_existence}. Finally, in Section \ref{Section_SPsystem} we derive from it the existence result for the Schr\"odinger-Poisson system \eqref{SP:system}.

\vskip0.2truecm
\paragraph{\textbf{Notation}} For $R>0$ and $x_0\in\R^N$ we denote by $B_R(x_0)$ the ball of radius $R$ and centre $x_0$, and we omit the centre when $x_0=0$. Given $\Omega\subset\R^N$, we denote $\Omega^c:=\R^N\setminus\Omega$, and its characteristic function by $\chi_\Omega$. The space of infinitely differentiable functions which are compactly supported is denoted by $C^\infty_c(\R^N)$, $\mathcal M(\R^N)$ stands for the space of measurable functions in $\R^N$, while $\cS$ is the Schwartz space of rapidly decreasing functions and $\cS'$ the dual space of tempered distributions. For $p\in[1,+\infty]$ the Lebesgue space of $p$-integrable functions is denoted by $L^p(\R^N)$ with norm $\|\cdot\|_p$. For $q>1$ we define its conjugate H\"older exponent as $q':=\frac q{q-1}$. The symbol $\lesssim$ indicates that an inequality holds up to a multiplicative constant depending only on structural constants, while $f \asymp g$ that $c_1f\leq g\leq c_2f$ for some $c_1,c_2>0$. Finally, $\text{o}_n(1)$ denotes a vanishing real sequence as $n\to+\infty$. Hereafter, the letter $C$ will be used to denote positive constants which are independent of relevant quantities and whose value may change from line to line.

\section{Functional space, functional inequalities, and preliminaries}\label{Section_spaces}

\subsection{The space $D^{1,2}L^p_\omega(\R^2)$}

Let us define the linear space
\begin{equation*}
	D^{1,2}(\R^2):= \left\{u \in L^1_\loc(\R^2)\,|\,\nabla u\in [L^2(\R^2)]^2\right\}
\end{equation*}
with seminorm $\|\nabla\cdot\|_2$. Note that, by the unboundedness of the domain, this seminorm cannot control the $L^2$-norm of the elements of $D^{1,2}(\R^2)$, and therefore $D^{1,2}(\R^2)\supsetneq H^1(\R^2)$. To retrieve a normed space, one needs to introduce the relation $u\sim v$ $\Leftrightarrow$ $v=u+c$ with $c\in\R$, and define the quotient space $\dot{D}^{1,2}(\R^2):=\left\{[u]_\sim\,|\,u\in D^{1,2}\right\}$, which turns out to be a Hilbert space with norm $\|\nabla\cdot\|_2$ (see \cite[Lemma II.6.2]{Galdi2011}). On the other hand, one may also introduce the space
\begin{equation*}
	D^{1,2}_0(\R^2):=\text{completion of } C^\infty_c(\R^2) \text{ w.r.t. } \|\nabla\cdot\|_2\,.
\end{equation*}
By \cite[Theorem II.7.5]{Galdi2011} the spaces $D^{1,2}_0(\R^2)$ and $\dot{D}^{1,2}(\R^2)$ are isomorphic.
\vskip0.2truecm
In order to find a suitable variational framework for the system \eqref{SP:system}, for $p>2$ and $b>1$ we define
$$
L^p_\omega(\R^2):=L^p(\R^2,\omega\dd x):=\{u\in\mathcal M(\R^2)\,|\,\|u\|_{\ast,p}<+\infty\},
$$
where the weight function $\omega$ is given by $\omega(x):=\ln(b+|x|)$, and
\begin{equation*}
	\|u\|_{\ast,p}:=\left(\int_{\R^2}|u|^p\ln(b+|x|)\dd x \right)^\frac1p,
\end{equation*}
and we consider the space
\begin{equation*}
	D^{1,2}L^p_\omega(\R^2):=D^{1,2}(\R^2) \cap L^p(\R^2,\omega\dd x)\,,
\end{equation*}
with norm 
\begin{equation*}
	\|u\|
	:=\left[\int_{\R^2}|\nabla u|^2 \dd x  + \left(\int_{\R^2}|u|^p\ln(b+|x|) \dd x \right)^{\frac{2}{p}}\right]^\frac12.
\end{equation*}
Note that, by the choice of $b>1$, in $D^{1,2}L^p_\omega(\R^2)$ it is possible to control the $L^p-$norm by the seminorm $\|\cdot\|_{\ast,p}$. Indeed,
\begin{equation}\label{p:norm:control:by:seminorm}
	\|u\|^p_p \leq(\ln b)^{-1}\int_{\R^2}\ln(b+|x|)|u(x)|^p \dd x =(\ln b)^{-1}\|u\|^p_{*,p}\,.
\end{equation}
Therefore,
\[
D^{1,2}(\R^2) \cap L^p(\R^2,\omega \dd x) = D^{1,2}(\R^2) \cap L^p(\R^2) \cap L^p(\R^2,\omega \dd x),
\]
and so, since $C^{\infty}_c(\R^2)$ is dense in both $D^{1,2}_0(\R^2)$ and $L^p(\R^2)$, we have the characterisation
\begin{align*} 
	D^{1,2}(\R^2)&\cap L^p(\R^2,\omega\dd x)= D^{1,2}(\R^2)\cap L^p(\R^2) \cap L^p(\R^2,\omega\dd x) \\
	&= D^{1,2}_0(\R^2) \cap L^p(\R^2) \cap L^p(\R^2,\omega\dd x)= D^{1,2}_0(\R^2) \cap L^p(\R^2,\omega\dd x).
\end{align*}
Furthermore, by \cite[Theorem 1.11]{KufnerOpic1984} $D^{1,2}L^p_\omega(\R^2)$ is a reflexive Banach space, whose dual is given by 
\begin{equation*}
	\begin{aligned}
		D^{-1,2}L^p_\omega(\R^2):&=(D^{1,2}(\R^2) \cap L^p_\omega(\R^2))'\\
		&=D^{-1,2}(\R^2)\big|_{D^{-1,2}L^p_\omega} + L^{p'}(\R^2,\omega\dd x)\big|_{D^{-1,2}L^p_\omega},
	\end{aligned}
\end{equation*}
since $\left(L^p_\omega(\R^2)\right)'=L^{p'}(\R^2,\omega\dd x)$, see \cite[Theorem 14.9]{Simon2017}, and also \cite[Theorem II.8.1]{Galdi2011} for the representation of the space $D^{-1,2}(\R^2)$. Let us now state important embedding properties of our space.
\begin{Prop}\label{cpt_emb_X}
	The space $X:=D^{1,2}L^p_\omega(\R^2)$ is compactly embedded in $L^q(\R^2)$ for all $q\geq p$.
\end{Prop}
\begin{proof}
	The embedding $X\hookrightarrow L^p(\R^2)$ is a consequence of \eqref{p:norm:control:by:seminorm}. Let now $q>p$, then by the Gagliardo-Nirenberg inequality (see Proposition \ref{GN} below, applied with $j=0, m=1, r=2, N=2, q=p$) we have
	\begin{equation*}
		\|u\|_q \leq C\|u\|^{\frac pq}_p\|\nabla u\|^{\frac{q-p}{q}}_2\,,
	\end{equation*}
	which implies $X\hookrightarrow L^q(\R^2)$, for every $q>p$.\\
	Let us now prove the compactness of these embeddings by relying on the Riesz criterion, see \cite[Theorem XIII.66]{ReedSimonIV1978}, which requires the continuity of the translation in the $L^q-$norm and a uniform decay at infinity of the elements in $X$. To this aim, let $S\subset X$ be a bounded subset, which is then also bounded in $L^q(\R^2)$ for $q\geq p$. Let also $R>0$ and $u\in S$. Then, by the H\"older inequality,
	\begin{equation*}
		\int_{\{|x|\geq R\}}|u|^q\dd x\leq\|u\|_{(q-1)p'}^{q-1}\left(\int_{\{|x|\geq R\}}|u|^p\dd x\right)^\frac1p\leq C\left(\int_{\{|x|\geq R\}}|u|^p\dd x\right)^\frac1p,
	\end{equation*}
	since $(q-1)p'\geq p$ and the continuity of the embedding shown before. Moreover,
	\begin{equation*}
		\int_{\{|x|\geq R\}}|u|^p\dd x\leq\frac{\|u\|_{\ast,p}^p}{\ln(b+R)}\leq\frac C{\ln(b+R)}\,,
	\end{equation*}
	since $u\in S$ bounded in $X$. Hence, for any $\varepsilon>0$ one can choose $R>0$ large enough such that
	\begin{equation}\label{Riesz_1}
		\int_{B_R^c}|u|^q\dd x\leq\varepsilon^q.
	\end{equation}
	Let us now prove the continuity of the translation in $L^q(\R^2)$. Since $X\subset D^{1,2}_0(\R^2)$, by density we can work within $C^\infty_c(\R^2)$. Let $u\in C^\infty_c(\R^2)$ and $h\in\R^2$. Following \cite[Proposition 9.3]{Brezis2011} and defining $\tau_hu:=u(\cdot+h)$, by the Jensen inequality we have
	\[
	|u(x+h)-u(x)|^2=\left|\int_0^1 h \cdot \nabla u(x+th) \dd t \right|^2 \leq |h|^2\int_0^1|\nabla u(x+th)|^2 \dd t.
	\]
	Integrating on $\R^2$ and using the Fubini-Tonelli theorem,
	\begin{equation*}
		\begin{aligned}
			\|\tau_hu - u\|_2^2&=\int_{\R^2}|u(x+h)-u(x)|^2 \dd x \\
			&\leq |h|^2\int_0^1\int_{\R^2}|\nabla u(x+th)|^2 \dd x \dd t = |h|^2\|\nabla u\|^2_2\,.
		\end{aligned}
	\end{equation*} 
	Hence, 
	\begin{align*}
		\|\tau_h u - u\|^q_q 
		&= \int_{\R^2}|u(x+h)-u(x)|^q \dd x\\
		&\leq \left(\int_{\R^2}|u(x+h)-u(x)|^{2(q-1)} \dd x\right)^\frac12\|\tau_hu-u\|_2\\
		&\les \left(\|\tau_hu\|_{2(q-1)}^{q-1} + \|u\|_{2(q-1)}^{q-1}\right)|h|\|\nabla u\|_2\\
		&\les 2\|u\|_{2(q-1)}^{q-1}\|u\||h|\les\|u\|^q|h|
	\end{align*}
	by the continuous embedding showed above, since $p>2$ implies that $2(q-1)>p$ for all $q\geq p$.\\
	The above inequality, together with \eqref{Riesz_1}, completes the proof.
\end{proof}

\subsection{Functional inequalities in log-weighted Lebesgue spaces}

Since our nonlinearities exhibit exponential growth, we present here some important functional inequalities, including newly established ones, to effectively handle them. The first result is a generalised Cao's inequality in $D^{1,2}(\R^2) \cap L^q(\R^2)$ taken from \cite{GuoLiu2023}, which is a particular case of the very general result obtained therein.
\begin{Th}[\cite{GuoLiu2023}, Theorem 1.1]\label{GuoLiu_Thm}
	Let $q \geq p$ and $\alpha < 4\pi$. Then there exists a constant $C:=C(p,q,\alpha) > 0$ such that for all $u \in D^{1,2}(\R^2) \cap L^p(\R^2)$ with $\|\nabla u\|_2 \leq 1$ there holds
	\begin{equation}\label{Cao:generalized}
		\int_{\R^2}|u|^q\,\e^{\alpha u^2}\dd x \leq C\|u\|^p_p\,.
	\end{equation}
	If $\alpha\geq 4\pi$, \eqref{Cao:generalized} remains true but the constant $C$ is not uniform in $u$.
\end{Th}

Next, we recall a Poho\v zaev-Trudinger inequality with logarithmic weight in $D^{1,2}L^p_\omega(\R^2)$ from \cite{CassaniTarsi2021}.
\begin{Th}[\cite{CassaniTarsi2021}, Theorem 3.3]\label{thm_wTp}
	Let $g$ satisfy assumptions $(f_1)$ and either $(f_{sc})$ or $(f_c)$ with exponent $\alpha_0$, and $G(s)=\int_0^s g(t) \dd t$. Then, the space $D^{1,2}L^p_{\omega}(\R^2)$ embeds into the weighted Orlicz space $L_G(\R^2,\omega\dd x)$, namely
	\begin{equation}\label{ineq_CaTa}
		\intd G(\alpha |u|)\ln(b+|x|)\dd x<+\infty
	\end{equation}
	for any $u\in D^{1,2}L^p_{\omega}(\R^2)$ and any $\alpha>0$.
\end{Th}
Reasoning as in \cite[Corollary 3.4]{CassaniTarsi2021}, from Theorem \ref{thm_wTp} it is easy to obtain the following:
\begin{Cor}\label{CaTa:Cor_3.4}
	For any $\alpha > 0$ the functional
	$$u \mapsto \int_{\R^2} G(\alpha|u|)\ln(b+|x|) \dd x\,, \quad u \in D^{1,2}L^p_\omega(\R^2)\,,$$
	is continuous, where $G$ is as in Theorem \ref{thm_wTp}.
\end{Cor}

For our purposes, a uniform estimate in the spirit of Moser for the inequality \eqref{ineq_CaTa} will be essential. In \cite[Theorem 3.3]{CassaniTarsi2021} the authors obtained it under the condition $\alpha\leq\frac{2\pi}{\alpha_0\sqrt p}$ by relying on Ruf's inequality. Here, we need to improve this threshold, in particular avoiding the dependence on $p$. In light of assumption ($f_c$), this is accomplished once the following Moser-type result is achieved, which extends to the space $D^{1,2}L^p_\omega(\R^2)$ the corresponding result in \cite{Tarsi} for the case $p=2$.
\begin{Th}\label{prop4pinonsh}
	For any $\alpha<4\pi$
	\begin{equation*}
		\sup_{\|u\|^2\leq 1}\intd|u|^p\e^{\alpha u^2}\ln(b+|x|)\dd x<+\infty\,.
	\end{equation*}
	The inequality is sharp, in the sense that, for any  $\alpha>4\pi$
	\begin{equation*}
		\sup_{\|u\|^2\leq 1}\intd|u|^p\e^{\alpha u^2}\ln(b+|x|)\dd x=+\infty\,.
	\end{equation*}
\end{Th}

\begin{proof}
	The proof follows by suitably adapting the arguments in \cite{Tarsi}. 
	The main tool is a a transformation which relates functions in the weighted space $D^{1,2}L^p_w(\R^2)$ to functions in the unweighted space $D^{1,2}_0(\R^2)\cap L^p(\R^2)$, based upon a change of variables acting only on the radial part of $x$. The price to pay is a dilation term in the Dirichlet norm, whose effect is lower and lower as $|x|\to +\infty$ due to the logarithmic growth of the weight. This property is the key tool that allows to retain the sharp threshold $4\pi$, and it is the only one which needs a proper delicate modification when passing from \cite{Tarsi} to our case. Let us denote by $T:\R^2\to\R^2$ the change of variable
	$$
	s=T(r)=\sqrt{2\int_0^r \rho\ln(b+\rho)\dd\rho}\,. 
	$$
	Setting
	\[
	v(y):=u(x), \   \ \ 	 \ \ \hbox{namely, } \ \ v( y)=u\left(T^{-1}(|y|)\cos \theta, T^{-1}(|y|)\sin \theta\right),
	\]
	as in \cite{Tarsi} one can easily verify that the map
	\begin{eqnarray*}
		\mathcal{T}: D^{1,2}L^p_\omega(\R^2)&\to& D^{1,2}_0(\R^2)\cap L^p(\R^2)\\
		u&\mapsto& v
	\end{eqnarray*}
	is an invertible, continuous map, with  continuous inverse map. Let $\chi(|x|)$ be a smooth, radial cut-off function
	$$\chi(|x|)=
	\begin{cases}
		1 &\hbox{if } |x|<1\,,\\
		0 &\hbox{if } |x|>2\,,
	\end{cases}
	$$
	and $\xi(|x|):=1-\chi(|x|)$; then scale $\xi$ as follows
	$$\xi_{\eta}(|x|):= \xi(\eta|x|)\quad\mbox{for}\ \eta\in(0,1)\,,$$
	and define $u_\eta:=u\cdot \xi_\eta$ for any $u\in X$ with $\|u\|\leq 1$; $u_\eta$ is then supported in $B_{1/\eta}^{\,c}$. We have
	\begin{equation}\label{dir_norm_ueta}
		\begin{aligned}
			\int_{\R^2}|\nabla u_\eta|^2
			&\leq \|\nabla u\|_2^2+\left(\int_{u^2\leq \eta^{\frac 2{p-2}}}|\nabla \xi_\eta|^2u^2+\int_{u^2>\eta^{\frac 2{p-2}}}|\nabla \xi_\eta|^2u^2\right)\\
			&\quad+2\left(\int_{u^2\leq\eta^{\frac 2{p-2}}}u\,\xi_\eta \nabla u \nabla \xi_\eta +\int_{u^2>\eta^{\frac 2{p-2}}}u\,\xi_\eta \nabla u \nabla \xi_\eta \right)\\
			&\leq \|\nabla u\|_2^2+\eta^{\frac 2{p-2}} \|\nabla \chi\|_2^2+\eta ^{-1}\int_{u^2>\eta^{\frac 2{p-2}}}|\nabla \xi_\eta|^2|u|^p\\
			&\quad +2\left(\eta^{\frac 1{p-2}}\|\nabla u\|_2 \|\xi_\eta \nabla \xi_\eta\|_2+\left[\int_{u^2>\eta^{\frac 2{p-2}}}|u|^2\right]^{\frac 12}\|\nabla u\|_2\| \nabla \xi_\eta\|_\infty\right)\\
			&\leq \|\nabla u\|_2^2+\eta^{\frac 2{p-2}} \|\nabla \chi\|_2^2+\eta \|\nabla \chi\|_\infty^2\int_{u^2>\eta^{\frac 2{p-2}}}|u|^p\\
			&\quad +2\left(\eta^{\frac 1{p-2}}\| \nabla u\|_2 \| \nabla \chi\|_2 +\sqrt\eta\|u\|_{\ast,p}^{\frac p2}\|\nabla u\|_2\| \nabla \chi\|_\infty\right)
		\end{aligned}
	\end{equation}
	so that, if $\tau:=\min\big\{\frac 12,\frac 1{p-2}\big\}$,
	\begin{equation}\label{norm_inex}
		\|u_\eta\|^2 \leq  1+C\eta^\tau
	\end{equation}
	for some positive fixed constant $C$ independent of $\eta$ and for any $u\in X$ such that $ \|u\|\leq 1$. The proof now follows as in \cite{Tarsi}. Let $v_\eta$ be defined by
	\[ 
	\ u_\eta(r\cos \theta, r\sin \theta)=v_\eta\left(T(r)\cos \theta, T(r)\sin \theta\right); 
	\]
	observing that 
	\begin{equation*}
		\|\nabla v_\eta\|_2^2+\|v_\eta\|_p^{\frac 2p}<\left(1+\frac 1{2\ln(1+\frac 1\eta)}\right)\|u_\eta\|^2\leq \left(1+\frac 1{2\ln(1+\frac 1\eta)}\right)(1+C\eta^\tau)
	\end{equation*}
	by  \eqref{norm_inex} we can  fix  $\eta$ such that
	$$
	\left(1+\frac 1{2\ln(1+\frac 1\eta)}\right)(1+C\eta^\tau)<\frac{4\pi}\alpha\,,
	$$
	so that $\alpha\|v_\eta\|^2<4\pi$, and apply Theorem \ref{GuoLiu_Thm} to conclude.\\
	The sharpness can be proved considering the same sequence of radial functions introduced in \cite{Tarsi} 
	\begin{equation*}
		u_n(x)=\frac1{\sqrt{2\pi}}\left\{
		\begin{array}{ll}
			\displaystyle \sqrt{\delta_n\ln n} & \displaystyle 0\leq |x|\leq\frac{R_n}{n}\,,\\
			\\
			\displaystyle \frac{\sqrt{\delta_n}}{\sqrt{\ln n}}\ln\left(\frac{R_n}{|x|}\right)
			& \displaystyle \frac{R_n}{n}<|x|\leq R_n\,, 
		\end{array}
		\right.
	\end{equation*}
	where
	\begin{equation*}
		R_n:=\frac{\sqrt{\ln n}}{\ln \ln n}\to+\infty\,, \quad\delta_n:= 1-\frac{\ln\ln n}{4\ln n}\to 1^-,\quad
		\text{ as } n\to+\infty\,.
	\end{equation*}
	Then 
	$$
	\|\nabla u_n\|_2^2=\delta_n\,,
	$$
	whereas, recursively integrating by parts,
	\begin{equation*}
		\begin{aligned}
			\|u_n\|_{\ast,p}^p&= \delta^{\frac p2}_n \ln^{\frac p2} n\int_0^{\frac{R_n}n}\!\!r\ln(b+r)\dd r+ \frac{\delta^{\frac p2}_n}{\ln^{\frac p2} n}\int_{\frac{R_n}n}^{R_n}\! r\ln^p\!\left(\frac{R_n}{r}\right)\ln(b+r)\dd r\\
			& \leq \frac{\ln^{\frac p2} n}{2n^2}R_n^2\ln\left(b+\frac{R_n}{n}\right)+\frac1{\ln^{\frac p2}n}\int_{\frac{R_n}n}^{\frac{R_n}\e}\!r\ln^p\left(\frac{R_n}{r}\right)\ln(b+r)\dd r\\
			&\quad +\frac{1}{\ln^{\frac p2} n}\int_{\frac{R_n}\e}^{R_n}r\ln(b+r)\dd r\\
			&\leq \frac{\ln^{\frac p2} n}{n^2}R_n^2 + C_p\bigg(\frac{R_n^2}{\ln^{\frac p2}n}\ln\left(b+\frac{R_n}\e\right)\\
			&\quad+\frac1{\ln^{\frac p2} n}\int_{\frac{R_n}n}^{\frac{R_n}\e}\!r\ln^{p-[p]-1}\!\left(\frac{R_n}{r}\right)\ln(b+r)\dd r\bigg) +\frac{R^2_n}{2\ln^{\frac p2} n}\ln\left(b+R_n\right)\\
			&\leq\frac{(\ln n)^{\frac{p+2}2}}{n^2\ln\ln n} + \frac{C_p\ln n}{2\ln^{\frac p2} n \cdot \ln \ln n} \to 0
		\end{aligned}
	\end{equation*}
	since $p>2$. On the other hand,
	\begin{equation*}
		\begin{aligned}
			\int_{\R^2}u_n^p \e^{4\pi u_n^2}\ln(b+|x|)\dd x&\geq 2\pi \int_0^{R_n/n}u_n^p\,\e^{2\delta_n\ln n}r\ln(b+r)\dd r \\
			&\geq\frac{\pi^{1-\frac p2}}{2^{\frac p2}}\ln b\,\delta_n^{\frac p2}\,\frac{(\ln n)^{\frac{p+1}2}}{(\ln\ln n)^2}\to +\infty \quad
			\text{ as } n\to+\infty\,.
		\end{aligned}
	\end{equation*}
\end{proof}

\begin{Rem}
	We remark that, although $\alpha=4\pi$ represents a sharp threshold for Theorem \ref{prop4pinonsh}, understanding the behaviour of the inequality at this critical value is far from trivial and requires a delicate analysis in the spirit of \cite{Tarsi}.
\end{Rem}

\noindent In view of Theorem \ref{prop4pinonsh}, the next two results are direct consequences.
\begin{Prop}\label{cpt_emb_X_log}
	The space $X:=D^{1,2}L^p_\omega(\R^2)$ is continuously embedded in $L^q_\omega(\R^2):=\{u\in\mathcal M(\R^2)\,|\,\intd |u|^q\ln(b+|x|)\dd x<+\infty\}$ for all $q\geq p$.
\end{Prop}
\begin{proof}
	The proof is a trivial consequence of Theorem \ref{prop4pinonsh} applied to any normalised function $u/\|u\| \in X$, combined with the elementary inequality
	$$
	|t|^{q}\leq C(q,p)|t|^p\e^{t^2}, \ \ t\in \R\,,
	$$
	which holds for any $q\geq p$ and for a proper constant $C(q,p)>1$.
\end{proof}

\begin{Th}
	Under the assumptions of Theorem \ref{thm_wTp}, for any $\alpha\leq\frac{4\pi}{\alpha_0}$ one has 
	\begin{equation*}
		\sup_{\|u\|^2\leq 1} \int_{\R^2} G\left(\alpha |u|\right)\ln(b+|x|)\dd x <+\infty\,.
	\end{equation*}
\end{Th}

\noindent Finally, we also prove a weighted Cao-type inequality in our setting.
\begin{Th}\label{propCaoweighted}
	For any $\alpha<4\pi$ and $M>0$
	\begin{equation*}
		\sup_{\|\nabla u\|_2\leq 1,\ \|u\|_{\ast, p}\leq M}
		\intd|u|^p\e^{\alpha u^2}\ln(b+|x|)\dd x=C(\alpha,M)<+\infty\,.
	\end{equation*}
\end{Th}
\begin{proof}
	The proof follows the same approach as the one for Theorem \ref{prop4pinonsh}, relying on the radial change of variable $T$ and on a quantitive estimate of the Dirichlet norm of the function $u_\eta$. As in \eqref{dir_norm_ueta}
	\begin{equation*}
		\begin{aligned}
			\int_{\R^2}|\nabla u_\eta|^2&\leq \|\nabla u\|_2^2+ \eta^{\frac 2{p-2}} \|\nabla \chi\|_{2}^2+\eta \|\nabla \chi\|^2_{\infty} M^p \\
			&\quad +2\left(\eta^{\frac 1{p-2}}\| \nabla u\|_2 \| \nabla \chi\|_2 +\sqrt\eta M^{\frac p2}\|\nabla u\|_2\| \nabla \chi\|_\infty\right)\\
			&\leq 1+ \eta^{\frac 2{p-2}} \|\nabla \chi\|_2^2+\eta \|\nabla \chi\|^2_\infty M^p +2\eta^{\frac1{p-2}}\|\nabla \chi\|_2 +2\sqrt\eta M^{\frac p2}\|\nabla \chi\|_\infty\\
			&\leq 1+C\eta^\sigma, \quad \sigma=\min\left\{\frac 1{p-2}, \frac 14\right\}
		\end{aligned}
	\end{equation*}
	if $\eta < M^{-2p}$. Let us then apply the change of variable $T$, obtaining a new function $v_\eta$; we easily obtain
	\begin{equation*}
		\|\nabla v_\eta\|_2^2<\left(1+\frac 1{2\ln(1+\frac 1\eta)}\right)\|\nabla u_\eta\|_2^2< 1+\frac C{|\ln\eta|}\,,
	\end{equation*}
	and
	$$
	\|v_\eta\|_p^2\leq \left(1+\frac C{|\ln \eta|}\right)\|u_\eta\|^2_{\ast, p}
	\leq \left(1+\frac C{|\ln \eta|}\right)\|u\|^2_{\ast,p}\leq \left(1+\frac C{|\ln \eta|}\right)M^2.
	$$
	Then, for any $\eta$ small enough, depending only on $M$ and $\alpha$, we get
	\begin{equation*}
		\begin{aligned}
			&\int_{B^{\,c}_{2/\eta}}|u|^p\e^{\alpha u^2}\ln(b+|x|)\dd x\leq \int_{\R^2}|u_\eta|^p\e^{\alpha  u_\eta^2}\ln(b+|x|)\dd x\\
			&\quad =\int_{\R^2}|v_\eta|^p\e^{\alpha  v_\eta^2}\dd x
			= \|\nabla v_\eta\|_2^p\int_{\R^2}\left(\frac{|v_\eta|}{\|\nabla v_\eta\|_2}\right)^p\e^{ \alpha \|\nabla v_\eta\|_2^2 \left(\frac{v_\eta}{\|\nabla v_\eta\|_2}\right)^2}\dd x\\
			&\quad \leq C(\alpha, \eta, M)
		\end{aligned}
	\end{equation*}
	if $\alpha \|\nabla v_\eta\|_2^2\leq 4\pi$, that is, if $1+\frac C{|\ln \eta|}<\frac{4\pi}{\alpha} $, where $\eta$ is now fixed such that $|\ln \eta|>C\left({\frac{4\pi}{\alpha}-1}\right)^{-1}$ and $\eta<M^{-2p}$.
\end{proof}

\subsection{Useful theorems and inequalities}

Throughout the paper, we will make great use of the following well-known results: the Hardy-Littlewood-Sobolev inequality (see \cite[Theorem 4.3]{LiebLoss}), and the Gagliardo-Nirenberg inequality (see \cite[eq. (2.2)]{Nirenberg1959}).
\begin{Prop}[HLS inequality]\label{HLS}
	Let $N\geq1$, $\mu\in(0,N)$, and $s,r>1$ with $\tfrac1s+\tfrac\mu N+\tfrac1r=2$. There exists a constant $C=C(N,\mu,s,r)$ such that for all $f\in L^s(\R^N)$ and $h\in L^r(\R^N)$ one has
	$$\int_{\R^N}\left(\frac1{|\cdot|^\mu}\ast f\right)\!h\dd x\leq C\|f\|_s\|h\|_r\,.$$
\end{Prop}
\begin{Prop}[GN inequality]\label{GN}
	Let $N\in\N$ and $u\in L^p(\R^N)$ such that $\nabla u\in L^r(\R^N)$, where $q,r \in [1,+\infty]$. Then, there exists a constant $C:=C(N,p,r,\theta)>0$ such that
	\begin{equation*}
		\|u\|_q\leq C\|\nabla u\|_r^\theta\|u\|_p^{1-\theta}.
	\end{equation*}
	where $\theta$ satisfies $\frac1p=\theta\left(\frac1r-\frac1N\right)+\frac{1-\theta}q$.
\end{Prop}
We also recall an abstract result from \cite[Proposition 3.1]{DuWeth2017} (see also \cite[Theorem 2.8]{Willem1996} for the version with the Palais-Smale condition), which will be needed to prove the existence of bounded Cerami sequences.
\begin{Prop}\label{abstract:MP}
	Let $\tX$ be a Banach space, $M_0$ be a closed subspace of a metric space $M$, and $\Gamma_0 \subset C(M_0,\tX)$. Define
	\[
	\tGamma=\left\{\gamma \in C(M,\tX) : \gamma_{\rvert_{M_0}} \in \Gamma_0\right\}.
	\]
	If $\Psi \in C^1(\tX,\R)$ satisfies
	\[
	\infty > c:=\inf_{\gamma \in \tGamma}\sup_{u \in M}\Psi(\gamma(u)) > a:=\sup_{\gamma_0 \in \Gamma_0}\sup_{u \in M_0}\Psi(\gamma_0(u))\,,
	\]
	then, for every $\eps \in \left(0,\frac{c-a}{a}\right)$, $\delta > 0$, and $\gamma \in \tGamma$ with $\sup_{u \in M} \Psi(\gamma(u)) \leq c + \eps$, there exists $u \in \tX$ such that
	\begin{itemize}
		\item $c-2\eps \leq \Psi(u) \leq c+2\eps$,
		\item $\dist(u,\gamma(M)) \leq 2\delta$,
		\item $(1+\|u\|_{\tX})\|\Psi'(u)\|_{\tX'} \leq \frac{8\eps}{\delta}$.
	\end{itemize}
\end{Prop}

\subsection{Consequences of the assumptions}\label{Assumption:consequences}
To end this section, let us point out some immediate consequences of ($f_1$)-($f_5$), ($f_{sc}$) and ($f_c$), which will be of use in our analysis, together with some comments in this regard:
\begin{enumerate}
	\item[i)] by $(f_1)$ and $(f_{sc})$, for any $r,\alpha>0$ and $s_0>1$ there is $C>0$ such that
	\begin{equation}\label{estFsc}
		0\leq F(s)\leq C\cdot\left\{
		\begin{array}{ll}
			s^p &\mbox{for } s\leq s_0\,,\\
			s^r\e^{\alpha s^2} &\mbox{for } s>s_0\,
		\end{array}\right.
	\end{equation}
	while, if $(f_{sc})$ is replaced by $(f_{c})$, the upper bound \eqref{estFsc} holds for $\alpha>\alpha_0$ and, moreover,
	\begin{equation}\label{estFc}
		0\leq f(s)\leq C\cdot
		\begin{cases}
			s^{p-1} &\mbox{for } s\leq s_0\,,\\
			s^r\e^{\alpha s^2} &\mbox{for } s>s_0\,;
		\end{cases}
	\end{equation}
	\item[ii)] by $(f_1)$ and $(f_{sc})$ or $(f_{c})$, there is $C>0$ such that for any $s>0$
	\begin{equation}\label{estFbelow}
		F(s)\geq C s^p\,;
	\end{equation}
	\item[iii)] assumption ($f_2$) implies that $f$ is monotone increasing. Moreover,
	$$\frac{\dd}{\dd t}\frac{F(t)}{f(t)}=\frac{f^2(t)-F(t)f'(t)}{f^2(t)}\leq 1-\tau\,,$$
	from which one infers
	\begin{equation}\label{ARcond_delta}
		F(t)\leq(1-\tau)tf(t)\quad\ \mbox{for any}\ \ t\geq0\,.
	\end{equation}
	\item[iv)] ($f_4$) is related to the well-known de Figueiredo-Miyagaki-Ruf condition \cite{dFMR} and is crucial in order to estimate the mountain pass level and gain compactness, see Lemma \ref{MP_level}. We note here that such an assumption, which goes back to \cite{dOdSdMSe}, avoids the prescription of a global lower-bound on $F$ of the kind \eqref{estFbelow} but with $C$ large enough: the latter, indeed, is widely used in the literature  but is not of practical verification. A condition similar to ($f_5$) appears also e.g. in \cite{CassaniTarsi2021,BCT,ChenShuTangWen2022,Romani1}.
	\item[v)] Examples of nonlinearities which satisfy the assumptions of Theorem \ref{existence:nontrivial:solution} are $F(s)=s^p$ and $F(s)=s^p\e^{s^\alpha}$ with $p>2$ and $\alpha\in(0,2)$ (concerning ($f_{sc}$)), and $F(s)=s^p\chi_{\{s<1\}}(s)+s^q\e^{s^2}\chi_{\{s\geq1\}}(s)$ with $q>-2$ (concerning ($f_c$)).
\end{enumerate}

\section{The variational framework}\label{Section_Var}

Formally, we can associate to the logarithmic Choquard equation \eqref{SP:equation} the energy functional $\cI:D^{1,2}L^p_\omega(\R^2)\to\R$ given by
\begin{equation}
	\label{energy:functional}
	\cI(u):=\frac12\int_{\R^2}|\nabla u|^2 \dd x + \frac12\int_{\R^2}\int_{\R^2} \ln|x-y|F(u(x))F(u(y)) \dd x \dd y\,.
\end{equation}
The aim of this section is to show that $\cI$ is indeed well-defined and regular in the space $D^{1,2}L^p_\omega(\R^2)$ described in Section \ref{Section_spaces}. First, we state an identity which will play a crucial role throughout the paper:
\begin{equation}\label{logarithmic:equality}
	\ln|x-y|=\ln(b+|x-y|)-\ln\left(1+\frac{b}{|x-y|}\right).
\end{equation}
This splitting was first used by \cite{Stubbe2008} with $b=1$, and subsequently developed by \cite{CingolaniWeth2017}. In \cite{ChenRadulescuWen,ChenShuTangWen2022} it was applied with $b>1$, and this allows for the embedding $L^p(\R^2,\omega\dd x)\hookrightarrow L^p(\R^2)$ as shown in Section \ref{Section_spaces}.

According to \eqref{logarithmic:equality}, and following the approach of \cite{CingolaniWeth2017}, we define the bilinear forms
\begin{gather*}
	(u,v) \mapsto A_1(u,v):=\int_{\R^2}\int_{\R^2} \ln(b+|x-y|)\,u(x)v(y) \dd x \dd y\,,\\
	(u,v) \mapsto A_2(u,v):=\int_{\R^2}\int_{\R^2} \ln\left(1+\frac{b}{|x-y|}\right)u(x)v(y) \dd x \dd y\,,\\
	(u,v) \mapsto A_0(u,v):=A_1(u,v) - A_2(u,v)=\int_{\R^2}\int_{\R^2} \ln|x-y|u(x)v(y) \dd x \dd y\,.
\end{gather*}
Since $b>1$, one has 
\begin{equation}\label{log_sum}
	\begin{aligned}
		\ln(b+|x-y|) &\leq \ln(b+|x|+|y|) \leq \ln(b+b|x|+b|y|) \\
		&\leq \ln((b+|x|)(b+|y|) = \ln(b+|x|) + \ln(b+|y|)\,,
	\end{aligned}
\end{equation}
and we can therefore estimate the bilinear form $A_1$ by
\begin{equation}
	\begin{aligned}\label{A1:estimate}
		|A_1(u,v)| 
		&\leq \int_{\R^2} \ln(b+|x|)|u(x)| \dd x\int_{\R^2}|v(y)| \dd y \\
		& \quad + \int_{\R^2}|u(x)| \dd x\int_{\R^2} \ln(b+|y|)|v(y)| \dd y \\
		& \leq \|u\|_{*,1}\|v\|_1+\|u\|_1\|v\|_{*,1}
	\end{aligned}
\end{equation}
for every $u,v \in L^1(\R^2,\omega\dd x)$. Concerning $A_2$, since $\ln(b+r) \leq r$ for every $r \geq 0$ (with the strict inequality if $r>0$), then, by the Hardy-Littlewood-Sobolev inequality (Proposition \ref{HLS}), one has
\begin{equation}\label{A2:estimate}
	A_2(u,v) \leq b\int_{\R^2}\int_{\R^2} \frac{1}{|x-y|}u(x)v(y) \dd x \dd y \les\|u\|_\frac43\|v\|_\frac43
\end{equation}
for every $u,v \in L^\frac43(\R^2)$. For $F\in C(\R)$, we also define the following functionals:
\begin{gather*}
	u \mapsto I_1(u):=A_1(F(u),F(u))=\int_{\R^2}\int_{\R^2} \ln(b+|x-y|)F(u(x))F(u(y)) \dd x \dd y \,,\\
	u \mapsto I_2(u):=A_2(F(u),F(u))=\int_{\R^2}\int_{\R^2} \ln\left(1+\frac{b}{|x-y|}\right)F(u(x))F(u(y)) \dd x \dd y \,,\\
	u \mapsto I_0(u):=A_0(F(u),F(u))=\int_{\R^2}\int_{\R^2} \ln|x-y|F(u(x))F(u(y)) \dd x \dd y \,.
\end{gather*}
Note that, if $F(u)\geq0$, then $I_1(u)\geq0$ and $I_2(u)\geq0$. With this notation, the energy functional \eqref{energy:functional} can be rewritten as
\[
\cI(u) = \frac12\|\nabla u\|^2_2 + \frac12I_0(u),
\]
and we are going to prove that $\cI$ is well-defined in $X:=D^{1,2}L^p_\omega(\R^2)$, see Proposition \ref{functional:regularity}. For the rest of the paper we always use this notation to indicate our space. 
\vskip0.2truecm
Before going into the details of the proof, we prepose an extension of \cite[Lemma 2.1]{CingolaniWeth2017} to our framework. It will be crucial in order to transfer estimates from the bilinear form $A_1$ to the norm $\|\cdot\|_{\ast,p}$, since it will be mainly applied with $\varphi_n=F(u_n)$ in combination with the lower bound \eqref{estFbelow}.

\begin{Lem}\label{Lemma3.1tilde}
	Let $p>1$, $u\in L^p(\R^2)\setminus\{0\}$ and nonnegative sequences $\{\varphi_n\}_n\subset L^1(\R^2)$, and $\{u_n\}_n\subset L^p(\R^2)$ such that $u_n\to u$ a.e. in $\R^2$ as $n\to+\infty$. Let moreover $F\in C(\R)$ with $F(t)>0$ for $t>0$.
	\begin{enumerate}
		\item[(a)] If $A_1(F(u_n),\varphi_n)\leq C$ and $\|\varphi_n\|_1\leq C$ for all $n\in\N$, then there exist $n_0\in\N$ and $\bar C>0$ such that $\|\varphi_n\|_{\ast,1}\leq\bar C$ for all $n>n_0$.
		\item[(b)] If $A_1(F(u_n),\varphi_n)\to0$ and $\|\varphi_n\|_1\to0$ as $n\to+\infty$, then $\|\varphi_n\|_{\ast,1}\to0$ as $n\to+\infty$.
	\end{enumerate}
\end{Lem}
\begin{proof}
	Since $u_n\to u$ a.e. in $\R^2$ and $F$ is continuous, by Egorov's theorem there exist $n_0\in\N$, $R>0$, and $\delta>0$, and a measurable set $A\subset B_R$ with positive measure, such that $F(u_n(x))\geq\delta$ for all $n\geq n_0$ and $x\in A$. For $x\in A$ and $y\in\R^2\setminus B_{(1+b)R}$ we have
	\begin{equation*}
		b+|x-y|\geq b+|y|-|x|\geq b+\left(1-\frac1{b+1}\right)|y|=b\left(1+\frac1{b+1}|y|\right)\geq b(1+|y|)^\frac1{b+1}
	\end{equation*}
	by Bernoulli's inequality. Hence,
	\begin{equation*}
		\begin{aligned}
			A_1(F(u_n),\varphi_n)&\geq
			\int_{\R^2\setminus B_{(b+1)R}} \int_A \ln(b+|x-y|)F(u_n(x))\varphi_n(y)\dd x\dd y\\
			&\geq 
			\left(\int_A F(u_n(x))\dd x\right)\!\left(\int_{\R^2\setminus B_{(b+1)R}}\!\!\ln\left(b(1+|y|)^\frac1{b+1}\right)\varphi_n(y)\dd y\right) \\
			&\geq 
			\frac{\delta|A|}{b+1}\int_{\R^2\setminus B_{(b+1)R}}\ln\left(b^{b+1}(1+|y|)\right)\varphi_n(y)\dd y\\
			&\geq 
			\frac{\delta|A|}{b+1}\int_{\R^2\setminus B_{(b+1)R}}\ln\left(b+|y|\right)\varphi_n(y)\dd y\\
			&\geq 
			\frac{\delta|A|}{b+1}\big(\|\varphi_n\|_{\ast,1}-\ln(b+(b+1)R)\|\varphi_n\|_1\big),
		\end{aligned}
	\end{equation*}
	having used the fact that $b^{b+1}>b>1$. This yields both (a) and (b), since then
	\begin{equation*}
		\|\varphi_n\|_{\ast,1} \leq \frac{b+1}{\delta|A|C}A_1(F(u_n),\varphi_n)+\ln(b+(b+1)R)\|\varphi_n\|_1\,.
	\end{equation*}
\end{proof}

\subsection{Regularity of the functional $\cI$}
Now, we move our attention to the well-posedness and regularity of the functional in our space $X=D^{1,2}L^p_\omega(\R^2)$.
\begin{Prop}\label{functional:regularity}
	Under ($f_1$) and ($f_{sc}$) or ($f_c$), the functionals $I_1,\,I_2,\,I_0,$ and $\cI$ are well-defined and of class $C^1$ on $X$, and moreover
	\begin{equation*}
		\begin{aligned}
			\cI'(u)[v]&=\intd\nabla u\cdot\nabla v \dd x+ 2A_0(F(u),f(u)v)\\
			&=\intd\nabla u\cdot\nabla v\dd x+\intd\left(\ln|\cdot|\ast F(u)\right)f(u)v\dd x, \quad \text{ for all } u,v \in X.
		\end{aligned}
	\end{equation*}
\end{Prop}
\begin{proof} First, note that $\ln(b+|x|)\geq\ln b>0$, since $b>1$, implies that
	\begin{equation}\label{finiteness_F}
		\int_{\R^2} F(u(x)) \dd x\leq (\ln b)^{-1}\int_{\R^2}\ln(b+|x|)F(u(x))\dd x<+\infty
	\end{equation}
	by Theorem \ref{thm_wTp}. Hence, from \eqref{A1:estimate} it follows that
	\begin{equation*}
		I_1(u)\leq 2\|F(u)\|_{*,1}\|F(u)\|_1<+\infty\,.
	\end{equation*}
	On the other hand, combining \eqref{A2:estimate} with \eqref{estFsc}, where $\alpha>0$ or $\alpha>\alpha_0$ if ($f_{sc}$) or ($f_c$) is assumed, respectively, and $r>\frac34p>\frac32>1$, one has
	\begin{equation}\label{I2:estimate}
		\begin{aligned}
			I_2(u)&\leq b\|F(u)\|_{\frac43}^2\les\|u\|_{\frac 43 p}^{2p}+\left(\int_{\R^2}|u|^{\frac43r}\e^{\frac43\alpha|u|^2}\dd x\right)^\frac32\\
			&\les\|u\|_{\frac 43 p}^{2p}+\|\nabla u\|_2^{2r}\left(\int_{\R^2}\left(\frac u{\|\nabla u\|_2}\right)^{\frac43r}\e^{\frac43\alpha\|u\|^2\left(\frac u{\|\nabla u\|_2}\right)^2}\dd x\right)^\frac32\\
			&\les\|u\|^{\frac32 p}_p\|\nabla u\|^{\frac p2}_2+C(u)\|\nabla u\|_2^{2r-\frac32p}\|u\|_p^{\frac32p}<+\infty\,,
		\end{aligned}
	\end{equation}
	by Proposition \ref{GN} and Theorem \ref{GuoLiu_Thm}. Consequently, $\cI$ is well-defined in $X$.
	
	Let now $\{u_n\}_n \subset X$ be a sequence such that $u_n \to u$ in $X$, that is
	\begin{equation}\label{un_conv}
		\|u_n-u\|^2 = \|\nabla u_n - \nabla u\|^2_2 + \|u_n-u\|^2_{*,p} \to 0\quad\text{as } n \to +\infty\,.
	\end{equation}
	We have
	\begin{align*}
		&|I_1(u_n) - I_1(u)| \\
		&
		\quad =\left|\int_{\R^2}\int_{\R^2} \ln(b+|x-y|)\big(F(u_n(x))F(u_n(y)) - F(u(x))F(u(y))\big) \dd x \dd y\right|\\
		&\quad \leq \int_{\R^2} \ln(b+|x|)F(u_n(x)) \dd x \int_{\R^2}\left|F(u_n(y)) - F(u(y))\right| \dd y \\
		&\qquad 
		+\int_{\R^2} \ln(b+|y|)\left|F(u_n(y)) - F(u(y))\right| \dd y\int_{\R^2}F(u_n(x)) \dd x \\
		&\qquad 
		+\int_{\R^2} \ln(b+|x|)\left|F(u_n(x)) - F(u(x))\right| \dd x \int_{\R^2}F(u(y)) \dd y \\
		&\qquad 
		+\int_{\R^2} \ln(b+|y|)F(u(y)) \dd y\int_{\R^2}\left|F(u_n(x)) - F(u(x))\right| \dd x\,,
	\end{align*}
	and all four terms tend to $0$ as $n\to\infty$, since the functionals $\int_{\R^2}F(u)\dd x$, $\int_{\R^2}F(u)\ln(b+|x|)\dd x$ on $X$ are continuous thanks to Corollary \ref{CaTa:Cor_3.4}.
	
	For any $u\in X$, the first G\^ateaux derivative of $I_1$ at $u$ along $v \in X$ is given by
	\[
	I_1'(u)[v] =2\int_{\R^2}\int_{\R^2} \ln(b+|x-y|)F(u(x))f(u(y))v(y) \dd x \dd y\,.
	\]
	By \eqref{log_sum}, we have
	\begin{align*}
		\frac12\left|I_1'(u)[v]\right| 
		&\leq
		\int_{\R^2}\ln(b+|x|)F(u(x)) \dd x \int_{\R^2}f(u(y))|v(y)| \dd y\\
		&\quad+\int_{\R^2}\ln(b+|y|)f(u(y))|v(y)| \dd y \int_{\R^2} F(u(x)) \dd x\\
		&\leq 
		\|F(u)\|_{*,1}\|f(u)\|_{\frac p{p-1}}\|v\|_p+\|F(u)\|_1\|f(u)\|_{*,\frac p{p-1}}\|v\|_{*,p}< +\infty
	\end{align*}
	by Theorems \ref{GuoLiu_Thm} and \ref{thm_wTp}. Now, let again $\{u_n\}_n \subset X$ and $u\in X$ be as in \eqref{un_conv}. We have
	\begin{align*}
		&\frac12\left|I_1'(u_n)[v]-I_1'(u)[v]\right|\\
		&\quad \leq
		\int_{\R^2} \ln(b+|x|)F(u(x)) \dd x \int_{\R^2}\left|f(u_n(y)) - f(u(y))\right||v(y)| \dd y \\
		&\qquad +
		\int_{\R^2} F(u(x)) \dd x \int_{\R^2}\ln(b+|y|)\left|f(u_n(y)) - f(u(y))\right||v(y)| \dd y \\
		&\qquad + 
		\int_{\R^2} \ln(b+|x|)\left|F(u_n(x)) - F(u(x))\right| \dd x \int_{\R^2}f(u_n(y))|v(y)| \dd y \\
		&\qquad + 
		\int_{\R^2} \left|F(u_n(x)) - F(u(x))\right| \dd x \int_{\R^2}\ln(b+|y|)f(u_n(y))|v(y)| \dd y\,.
	\end{align*}
	Recall now \eqref{finiteness_F}, and analogously
	\begin{equation}\label{cont_F}
		\begin{aligned}
			&\int_{\R^2}\left|F(u_n(x)) - F(u(x))\right|\dd x \\
			&\qquad \leq(\ln b)^{-1}\int_{\R^2}\ln(b+|x|)\left|F(u_n(x)) - F(u(x))\right|\dd x=\text{o}_n(1)
		\end{aligned}
	\end{equation}
	by Corollary \ref{CaTa:Cor_3.4}. Moreover, by H\"older's inequality
	\begin{equation*}
		\begin{aligned}
			\int_{\R^2}f(u_n(y))|v(y)|\dd y&\leq C(b)\int_{\R^2}\ln(b+|y|)f(u_n(y))|v(y)|\dd y\\
			&\les\|f(u_n)\|_{*,\frac p{p-1}}\|v\|_{\ast,p}\leq C(u)\|v\|\,,
		\end{aligned}
	\end{equation*}
	and
	\begin{align*}
		\begin{aligned}
			\int_{\R^2}\Big|f(u_n(y))&-f(u(y))\Big||v(y)|\dd y\\
			&\leq(\ln b)^{-1}\int_{\R^2}\ln(b+|y|)\left|f(u_n(y)) - f(u(y))\right||v(y)| \dd y\\
			&\les\|f(u_n)-f(u)\|_{*,\frac p{p-1}}\|v\|_{\ast,p}=\text{o}_n(1)\|v\|\,,
		\end{aligned}
	\end{align*}
	since $f$ is continuous. Combining the above inequalities, one infers
	\[
	\frac12\left|I_1'(u_n)[v]-I_1'(u)[v]\right| \leq C(u)\|v\|\,\text{o}_n(1)\,,
	\]
	namely $I_1'\in C(X)$. Let us now focus on $I_2$. For $\{u_n\}_n \subset X$ and $u\in X$ for which \eqref{un_conv} holds, by \eqref{A2:estimate} one has 
	\begin{align*}
		&|I_2(u_n)-I_2(u)|\\
		&\quad \leq 
		\int_{\R^2}\int_{\R^2} \ln\left(1+\frac{b}{|x-y|}\right)F(u_n(x))|F(u_n(y))-F(u(y))| \dd x \dd y \\
		&\qquad 
		+ \int_{\R^2}\int_{\R^2} \ln\left(1+\frac{b}{|x-y|}\right)F(u(y))|F(u_n(x))-F(u(x))| \dd x \dd y \\
		&\quad \les\|F(u_n)-F(u)\|_{\frac43}\left(\|F(u_n)\|_{\frac43} + \|F(u)\|_{\frac43}\right),
	\end{align*}
	which tends to $0$ as $n \to +\infty$, since $\|F(u_n)-F(u)\|_{\frac43}=\text{o}_n(1)$ by \eqref{cont_F}, $\|F(u)\|_{\frac43}<+\infty$ by \eqref{finiteness_F}, and $\|F(u_n)\|_{\frac43}\leq C$ by continuity.\\
	Computing the first variation of $I_2$ at $u \in X$ along $v \in X$ we get
	\[
	I_2'(u)[v]=2\int_{\R^2}\int_{\R^2} \ln\left(1+\frac{b}{|x-y|}\right)F(u(x))f(u(y))v(y) \dd x \dd y
	\]
	and so
	\begin{align*}
		\frac12\left|I_2'(u)[v]\right|&\leq b\int_{\R^2}\int_{\R^2}\frac1{|x-y|}F(u(x))f(u(y))|v(y)| \dd x \dd y \\
		&\les\|F(u)\|_{\frac43}\|f(u)\|_{\frac43p'}\|v\|_p\leq\|F(u)\|_{\frac43}\|f(u)\|_{\frac43p'}\|v\|< +\infty
	\end{align*}
	by Theorem \ref{GuoLiu_Thm}. Analogously, for $\{u_n\}_n \subset X$ and $u\in X$ as in \eqref{un_conv},
	\begin{multline*}
		\frac12\left|I_2'(u_n)[v]-I_2'(u)[v]\right|\\
		\les\left(\|F(u_n)\|_\frac43\|f(u_n) - f(u)\|_{\frac43p'}+\|F(u_n)-F(u)\|_\frac43\|f(u)\|_{\frac43p'}\right)\|v\|\,,
	\end{multline*}
	which again tends to $0$ by the above arguments. As a result, both $I_1$ and $I_2$ are of class $C^1$ on $X$: consequently, also $I_0=I_1-I_2$ and $\cI$ have the same regularity.
\end{proof}

\section{Analysis of Cerami sequences}\label{Section_Cerami}
Usually, a mountain pass geometry of the functional (see Lemma \ref{finite:positive:mp:level} below) directly provides the existence of a Cerami sequence, namely a sequence $\{u_n\}_n\subset X$ such that
\begin{equation}\label{Cerami_def}
	\cI(u_n) \to c_{mp}, \ \quad\  (1+\|u_n\|)\|\cI'(u_n)\|_{X'} \to 0,
\end{equation}
which yields the existence of a weak solution, by using some compactness argument which exploits the boundedness of such a sequence in $X$. In our case, however, the proof of the boundedness of the Cerami sequence is not standard, and we need to improve the properties that such a sequence has. The abstract result contained in Proposition \ref{abstract:MP} allows us to find a Cerami sequence with the additional property that $\cJ(u_n)\to0$ as $n\to+\infty$ (see Lemma \ref{existence:Cerami:sequence} below), where the functional $\cJ:X \to \R$ is given by
\begin{equation}\label{J:functional}
	u \mapsto \cJ(u):=2\|\nabla u\|^2_2 - 2 I_0(u) + 2 A_0(F(u), f(u)u)- \frac12\|F(u)\|^2_1\,.
\end{equation}
The boundedness of $\{u_n\}_n$ in $X$ will follow then by combining \eqref{Cerami_def} and $\cJ(u_n)\to0$. This strategy, which is reminiscent of the construction of Palais-Smale-Poho\v zaev sequences in the context of problems with prescribed mass \cite{Jeanjean}, was first employed in \cite{Ruiz2006} in the higher-dimensional case, and implemented in the planar setting by \cite{DuWeth2017,ChenShuTangWen2022} in the case $f(u)=u$. Here we need to extend it to the case of a general nonlinearity.
\vskip0.2truecm
Let us first define
\[
\Gamma:=\left\{\gamma \in C\left([0,1],X\right) : \gamma(0)=0 \text{ and } \cI(\gamma(1)) < 0 \right\},
\]
and the mountain pass level
\begin{equation}
	\label{mountain:pass:level}
	c_{mp}:=\inf_{\gamma \in \Gamma}\max_{t \in [0,1]}\cI(\gamma(t))\,.
\end{equation}

\begin{Lem}\label{finite:positive:mp:level}
	Assume $(f_1)$ and $(f_{sc})$ or $(f_c)$ hold. Then, the set $\Gamma$ is nonempty and $0 < c_{mp} < +\infty$.
\end{Lem}
\begin{proof}
	We start with a control from below of $I_1$: since $\ln(b+r) \geq \ln b$, for $r \geq 0$ and $b>1$, it follows by \eqref{estFbelow} that
	\begin{equation}\label{bound:I1:below}
		\begin{aligned}
			I_1(u) &= \int_{\R^2} \ln(b+|x-y|)F(u(x))F(u(y)) \dd x \dd y\\
			&\geq \ln b \left(\int_{\R^2}F(u(x)) \dd x \right)^2 
			\geq C \ln b \|u\|^{2p}_p.
		\end{aligned}
	\end{equation}
	Concerning $I_2$, let us refine the upper bound \eqref{I2:estimate}. By \eqref{estFsc} one has
	\begin{equation*}
		\begin{aligned}
			I_2(u)&\les\|u\|_{\frac 43 p}^{2p}+\|\nabla u\|_2^{2r}\left(\int_{\R^2}\left(\frac u{\|\nabla u\|_2}\right)^{\frac43r}\e^{\frac43\alpha\|u\|^2\left(\frac u{\|\nabla u\|_2}\right)^2}\dd x\right)^\frac32\\
			&\les\|u\|^{\frac32 p}_p\|\nabla u\|^{\frac p2}_2+\|\nabla u\|_2^{2r-\frac32p}\|u\|_p^{\frac32p}\,,
		\end{aligned}
	\end{equation*}
	having used Proposition \ref{GN} on the first term and Theorem \ref{GuoLiu_Thm} for the second, with the choice $r>\frac34p>\frac32>1$, $\alpha>\alpha_0$ close to $\alpha_0$ (with a little abuse of notation, for $\alpha_0=0$ if ($f_{sc}$) is assumed) and having required that $\frac43\alpha\|u\|^2<4\pi$. As a result, for any $u\in X$ with $\|u\|<\sqrt{3\pi\alpha^{-1}}$, and choosing now $r=p$, we get 	
	\begin{equation}\label{I2:finer:estimate}
		I_2(u)\les\|u\|^{\frac32p}_p\|\nabla u\|^{\frac p2}_2<+\infty\,.
	\end{equation}
	Hence, combining \eqref{bound:I1:below} and \eqref{I2:finer:estimate}, by Young's inequality with $\nu$ and $\nu'$ to be chosen later, we get
	\begin{equation}\label{I:nonnegative}
		\begin{aligned}
			\cI(u) &= \frac12\|\nabla u\|^2_2 + I_1(u) -I_2(u)\\
			&\geq \left(\frac12\|\nabla u\|^2_2 - \frac C{\nu'}\|\nabla u\|^{\frac p2\nu'}_2\right) 
			+ \left(\ln b\|u\|^{2p}_p - \frac C{\nu}\|u\|^{\frac32 p\nu}_p \right).
		\end{aligned}
	\end{equation}
	Choosing $\nu\in\left(\frac43,\frac4{(4-p)_+}\right)$, which is nonempty since $p>1$, we easily infer from \eqref{I:nonnegative} that $0$ is a local minimum for $\cI$.
	
	Let us now evaluate the functional $\cI$ along the fiber set $\{t^2u(t\,\cdot) : u \in X, t>0\}$. For a fixed $u\in X$ we have
	\begin{align*}
		\cI(t^2u(tx)) & =\frac{t^4}2\int_{\R^2}|\nabla u(x)|^2 \dd x\\
		&\quad + \frac{t^{-4}}2\int_{\R^2}\int_{\R^2}\ln|x-y|F(t^2u(x))F(t^2u(y))\dd x\dd y\\
		&\quad-\frac{t^{-4}\ln t}2\left(\int_{\R^2}F(t^2u(x))\dd x\right)^2.
	\end{align*}
	For $u\in C^\infty_c(B_{1/4}(0))$, one then gets
	\begin{equation*}
		\begin{aligned}
			\cI(t^2u(tx))&\leq\frac{t^4}{2}\int_{\R^2}|\nabla u(x)|^2\dd x-\frac{t^{-4}}{2}\left(\ln2+\ln t\right)\left(\int_{\R^2}F(t^2u(x))\dd x\right)^2\\
			&\leq\frac{t^4}{2}\int_{\R^2}|\nabla u(x)|^2\dd x-\frac{\ln2+\ln t}2t^{4(p-1)}\left(\int_{\R^2}|u|^p\dd x\right)^2\to-\infty
		\end{aligned}
	\end{equation*}
	as $t\to+\infty$. Hence $\sup_{t>0}\cI(t^2u(tx))<+\infty$, and there exists $t_*=t_*(u)>0$ such that $\cI(t_*^2u(t_*x))=\max_{t>0}\cI(t^2u(tx))$. Now, the function $\gamma(t)=(\tilde{t}t)^2u(\tilde{t}t \cdot)$, with $\tilde{t} >> t_*$ has the properties that $\gamma \in C([0,1],X)$, $\gamma(0)=0$, and $\cI(\gamma(1))<0$. As a result, $\gamma\in\Gamma$, namely $\Gamma\neq\emptyset$ and $c_{mp} < +\infty$.\\
	Since $\cI$ has a local minimum in $0$ by \eqref{I:nonnegative}, there exist a constant $a_0>0$ and $\rho>0$ such that 
	\[
	\cI(u) \geq a_0 \ \text{ if }\  u \in S_\rho(0) := \left\{u \in X : \|\nabla u\|^2_2 + \|u\|^p_p = \rho\right\}.
	\]
	Let $\gamma \in \Gamma$, then $\|\nabla\gamma(1)\|^2_2 + \|\gamma(1)\|^p_p > \rho$, and by the mean value theorem there exists $\bar{t} \in [0,1]$ such that $\|\nabla\gamma(\bar{t})\|^2_2 + \|\gamma(\bar{t})\|^p_p = \rho$. This means that $\gamma(\bar{t}) \neq 0$, hence $\cI(\gamma(\bar{t})) \geq a_0$. Therefore,
	\[
	\sup_{t \in [0,1]}\cI(\gamma(t)) \geq \cI(\gamma(\bar{t})) \geq a_0 > 0\,.
	\]
	Taking the infimum on $\Gamma$, we can conclude that $c_{mp}>0$.
\end{proof}

With the help of the abstract result of Proposition \ref{abstract:MP}, we are in a position to prove the existence of a specific Cerami sequence.
\begin{Lem}
	\label{existence:Cerami:sequence}
	Assume $(f_1)$ and $(f_{sc})$ or $(f_c)$ hold. Then there exists a Cerami sequence $\{u_n\}_n \subset X$ at the mountain pass level $c_{mp}$ defined in \eqref{mountain:pass:level}, such that 
	\begin{equation}\label{Jto0}
		\cJ(u_n)\to0\,.
	\end{equation}
\end{Lem}
\begin{proof}
	Let $\tilde{X}:=\R \times X$ be the Banach space endowed with the norm $\|(s,v)\|_{\tilde{X}}:=\left(|s|^2+\|v\|^2\right)^\frac12$. Consider the continuous map $\rho:\tilde{X} \to X$ defined as
	\begin{equation*}
		\rho(s,v)[x]:=\e^{2s}v(\e^sx)\,,\quad s \in \R,\ v \in X,\ x \in \R^2
	\end{equation*}
	and
	\[
	\Psi:=\cI \circ \rho : \tilde{X} \to \R\,.
	\]
	We compute
	\begin{align*}
		\Psi(s,v)&=\cI(\rho(s,v))\\
		&=\frac12\int_{\R^2}|\nabla \rho(s,v)|^2 \dd x  \\
		&\quad + \frac12\int_{\R^2}\int_{\R^2} \ln|x-y|F(\rho(s,v)[x])F(\rho(s,v)[y])\dd x \dd y \\
		&=
		\frac{\e^{6s}}{2}\int_{\R^2} |\nabla v(\e^sx)|^2 \dd x \\
		&\quad + \frac12\int_{\R^2}\int_{\R^2} \ln|x-y|F(\e^{2s}v(\e^sx))F(\e^{2s}v(\e^sy)) \dd x \dd y \\
		&=\frac{\e^{4s}}{2}\int_{\R^2} |\nabla v|^2 \dd x' \\
		&\quad + \frac{\e^{-4s}}{2}\int_{\R^2}\int_{\R^2} \ln|x'-y'|F(\e^{2s}v(x'))F(\e^{2s}v(y')) \dd x' \dd y' \\
		& \quad - s\frac{\e^{-4s}}{2}\left(\int_{\R^2}F(\e^{2s}v(x')) \dd x' \right)^2.
	\end{align*}
	By Lemma \ref{functional:regularity}, $\Psi$ is of class $C^1$ on $\tilde{X}$, therefore we can compute the partial derivatives of $\Psi$. For $s\in\R$ and $v\in X$ one has
	\begin{align*}
		\partial_s\Psi(s,v)
		&=2\e^{4s}\int_{\R^2}|\nabla v|^2\dd x-2\e^{-4s}\int_{\R^2}\int_{\R^2}\ln|x-y|F(\e^{2s}v(x))F(\e^{2s}v(y))\dd x\dd y \\
		&\quad+2\e^{-4s}\int_{\R^2}\int_{\R^2}\ln|x-y|F(\e^{2s}v(x))f(\e^{2s}v(y))\e^{2s}v(y)\dd x\dd y\\
		&\quad+2s\e^{-4s}\left(\int_{\R^2}F(\e^{2s}v)\right)^2-\frac{\e^{-4s}}2\left(\int_{\R^2}F(\e^{2s}v)\right)^2\\
		&\quad -2s\e^{-4s}\left(\int_{\R^2}F(\e^{2s}v)\right)\left(\int_{\R^2}f(\e^{2s}v)\e^{2s}v\right)\\
		&=2\e^{4s}\int_{\R^2}|\nabla v|^2\dd x-2\int_{\R^2}\int_{\R^2}\ln|x-y|F(\e^{2s}v(\e^sx))F(\e^{2s}v(\e^sy))\dd x\dd y\\
		&\quad+2\int_{\R^2}\int_{\R^2}\ln|x-y|F(\e^{2s}v(\e^sx))f(\e^{2s}v(\e^sy))\e^{2s}v(\e^sy)\dd x\dd y\\
		&\quad-\frac12\left(\int_{\R^2}F(\e^{2s}v(\e^sx))\dd x\right)^2\\
		&=\cJ(\rho(s,v))\,,
	\end{align*}
	where $\cJ$ is defined in \eqref{J:functional}. On the other hand, for $w\in X$, one has
	\begin{align*}
		\partial_v\Psi(s,v)[w] &= \partial_v\cI(\rho(s,v))[w]\\
		&= \e^{4s}\int_{\R^2} \nabla v(x) \cdot \nabla w(x) \dd x \\
		&\quad- s\e^{-4s}\left(\int_{\R^2}F(\e^{2s}v(x))\dd x\right)\left(\int_{\R^2}f(\e^{2s}v(x))w(x)\dd x\right)\\
		& \quad
		+ \e^{-4s}\int_{\R^2}\int_{\R^2} \ln|x-y|F(\e^{2s}v(x))f(\e^{2s}v(y))\e^{2s}w(y) \dd x \dd y \\
		&=\int_{\R^2} \nabla(\e^{2s}v(\e^sx')) \cdot \nabla(\e^{2s}w(\e^sx')) \dd x' \\
		& \quad
		- s\left(\int_{\R^2}F(\e^{2s}v(\e^sx'))\dd x'\right)\left(\int_{\R^2}f(\e^{2s}v(\e^sx'))w(\e^sx')\dd x'\right)\\
		& \quad
		+ \int_{\R^2}\int_{\R^2}\ln|x'-y'|F(\e^{2s}v(\e^sx'))f(\e^{2s}v(\e^sy'))\e^{2s}w(\e^sy')\dd x'\dd y'\\
		& \quad
		+ s\left(\int_{\R^2}F(\e^{2s}v(\e^sx'))\dd x'\right)\left(\int_{\R^2}f(\e^{2s}v(\e^sx'))w(\e^sx')\dd x'\right)\\
		&=\int_{\R^2} \nabla \rho(s,v) \cdot \nabla\rho(s,w) \dd x' \\
		& \quad
		+ \int_{\R^2}\int_{\R^2}\ln|x'-y'|F(\rho(s,v)[x'])f(\rho(s,v)[y'])\rho(s,w)[y']\dd x'\dd y'\\
		&=\cI'(\rho(s,v))[\rho(s,w)]\,.
	\end{align*} 
	Hence, the first variation of $\Psi$ at $(s,v) \in \tilde{X}$ along $(h,w) \in \tilde{X}$ is given by
	\begin{equation}
		\label{derivative:of:phi}
		\Psi'(s,v)(h,w) = \cI'(\rho(s,v))[\rho(s,w)] + \cJ(\rho(s,v))h\,.
	\end{equation}\\
	We are now going to apply Proposition \ref{abstract:MP} to the functional $\Psi$. To this end, let $\tilde{\Gamma}:=\{\tilde{\gamma} \in C([0,1],\tilde{X}) : \tilde{\gamma}(0)=(0,0),\ \Psi(\tilde{\gamma}(1))<0\}$ and
	\begin{equation}
		\label{phi:minmax:level}
		\tilde{c}:=\inf_{\tilde{\gamma} \in \tilde{\Gamma}}\max_{t \in [0,1]}\Psi(\tilde{\gamma}(t))\,.
	\end{equation}
	With these choices, it follows that
	\[
	(\rho \circ \tilde{\gamma})(0) = \rho(\tilde\gamma(0)) = \rho(0,0) = 0
	\]
	and
	\[
	\cI((\rho \circ \tilde{\gamma})(1)) = (\cI \circ \rho)(\tilde{\gamma}(1)) = \Psi(\tilde{\gamma}(1)) < 0\,,
	\]
	that is, $\Gamma=\{\rho \circ \tilde{\gamma} : \tilde{\gamma} \in \tilde{\Gamma}\}$ and the values \eqref{mountain:pass:level} and \eqref{phi:minmax:level} coincide. Let now $\{\gamma_n\}_n\subset\Gamma$ be a sequence of paths such that
	\[
	\sup_{t \in [0,1]} \cI(\gamma_n(t)) \leq c_{mp} + \frac{1}{n^2}\,.
	\]
	Defining $\widetilde\gamma_n(t):=(0,\gamma_n(t))$, which belongs to $\tilde{\Gamma}$, we have
	\[
	\sup_{t \in [0,1]}\Psi(\widetilde\gamma_n(t))= \sup_{t \in [0,1]} \cI(\gamma_n(t)) \leq c_{mp} + \frac1{n^2}\,.
	\]
	Hence, Proposition \ref{abstract:MP} applied with $M=[0,1]$ and $M_0=\{0,1\}$ yields the existence of a sequence $(s_n,v_n) \in \tilde{X}$ such that
	\begin{enumerate}[label=\textbf{(\alph*)}]
		\item $\Psi(s_n,v_n) \to \tilde c=c_{mp}$\,, \label{abstract:thesis:1}
		\item $\dist((s_n,v_n),\{0\} \times \gamma_n([0,1])) \to 0$\,, \label{abstract:thesis:2}
		\item $(1+\|(s_n,v_n)\|_{\tilde{X}})\|\Psi'(s_n,v_n)\|_{\tilde{X}'} \to 0$\,, \label{abstract:thesis:3}
	\end{enumerate}
	as $n \to +\infty$. We observe that \ref{abstract:thesis:2} implies 
	\begin{equation}
		\label{sn:to:zero}
		s_n \to 0\ \ \text{ as }\ \,n \to +\infty\,.
	\end{equation}
	Defining now $u_n:=\rho(s_n,v_n)$, by \ref{abstract:thesis:1} we get
	\[
	\cI(u_n) = \cI(\rho(s_n,v_n)) = \Psi(s_n,v_n) \to c_{mp}\,, \text{ as } n \to +\infty
	\]
	while, taking $h=1$ and $w=0$ in \eqref{derivative:of:phi}, from \ref{abstract:thesis:3} we also infer
	\begin{equation*}
		\cJ(u_n)=\cJ(\rho(s_n,v_n)) \to 0
	\end{equation*}
	as $n \to +\infty$. To obtain the last required property, observe that for a given $v \in X$, defining
	\begin{equation*}
		w_n=\e^{-2s_n}v(\e^{-s_n}\cdot),
	\end{equation*}
	\eqref{sn:to:zero} allows to show that
	\begin{equation}
		\label{wn:norm:asymptotics}
		\begin{aligned}
			&\|w_n\|^2 
			= \|\nabla w_n\|^2_2 
			+ \|w_n\|^2_{*,p} \\
			&= 
			\e^{-6s_n}\int_{\R^2} |\nabla v(\e^{-s_n}x)|^2 \dd x 
			+ \e^{-2s_n}\left(\int_{\R^2} \ln(b+|x|)|v(\e^{-s_n}x)|^p \dd x \right)^\frac2p\\
			&=\e^{-4s_n}\int_{\R^2} |\nabla v(x')|^2 \dd x'
			+ \e^{-2s_n\left(\frac{p-1}{p}\right)}\left(\int_{\R^2}\!\ln(b+\e^{s_n}|x'|)|v(x')|^p \dd x'\right)^\frac2p\\
			&=
			(1+\text{o}_n(1))\left(\int_{\R^2} |\nabla v(x')|^2 \dd x'+ \left(\int_{\R^2}\!\ln(b+\e^{s_n}|x'|)|v(x')|^p \dd x'\right)^\frac2p\right)\\
			&=
			(1+\text{o}_n(1))\|v\|^2 + \text{o}_n(1)\|v\|^2_p,
		\end{aligned}
	\end{equation}
	as $n\to+\infty$, where in the last step we used the fact that
	\begin{equation*}
		\begin{aligned}
			\ln(b+\e^{-s_n}|x|)&\leq \ln(b+|x|) + \ln\left(1+\frac{(1-\e^{-s_n})|x|}{b+|x|}\right)\\
			&\leq \ln(b+|x|) + \ln(1+(1-\e^{-s_n}))\,.
		\end{aligned}
	\end{equation*}
	Analogously one can show that $\|v_n\|=(1+\text{o}_n(1))\|u_n\|$, therefore, on the one hand, by \eqref{derivative:of:phi} with $h=0$ and \eqref{sn:to:zero}, one infers
	\begin{equation}\label{I'_1}
		\begin{aligned}
			(1+\|(s_n,v_n)&\|_{\tilde{X}})|\Psi'(s_n,v_n)(0,w_n)|\\
			&=(1+(|s_n|^2+\|v_n\|^2)^\frac12)|\cI'(\rho(s_n,v_n))\rho(s_n,w_n)|\\
			&=\left(1+\text{o}_n(1)+(1+\text{o}_n(1))\|u_n\|)\right)|\cI'(u_n)v|\,,
		\end{aligned}
	\end{equation}
	while, on the other hand, by \eqref{wn:norm:asymptotics},
	\begin{equation}\label{I'_2}
		\begin{aligned}
			(1+\|(s_n,v_n)\|_{\tilde{X}})|\Psi'(s_n,v_n)(0,w_n)|&\leq(1+\|(s_n,v_n)\|_{\tilde{X}})\|\Psi'(s_n,v_n)\|_{\tilde{X}'}\|w_n\|\\
			&=\text{o}_n(1)(1+\text{o}_n(1))\|v\|\,.
		\end{aligned}
	\end{equation}
	Combining together \eqref{I'_1}-\eqref{I'_2} we deduce
	$$
	\left(1+\text{o}_n(1)+(1+\text{o}_n(1))\|u_n\|\right)\|\cI'(u_n)\|_{X'}\to0
	$$
	as $n\to+\infty$, which readily implies
	$$
	\left(1+\|u_n\|\right)\|\cI'(u_n)\|_{X'}\to0
	$$
	as $n\to+\infty$.
\end{proof}
The extra property \eqref{Jto0} obtained in Proposition \ref{existence:Cerami:sequence} is crucial to prove the boundedness of a Cerami sequence in $X$, as shown next.
\begin{Lem}\label{boundedness:Cerami}
	Assume $(f_1)$, $(f_2)$ and $(f_{sc})$ or $(f_c)$ hold.
	Let $\{u_n\}_n \subset X$ be a sequence such that
	\begin{equation}\label{un_Cerami}
		\cI(u_n) \to c_{mp}\,, \quad (1+\|u_n\|)\|\cI'(u_n)\|_{X'} \to 0\,, \quad \cJ(u_n) \to 0\,.
	\end{equation}
	Then, $\{u_n\}_n$ is bounded in $X$. Moreover, there exist $C_1,C_2>0$ such that
	\begin{gather}
		|I_1(u_n)|=\left|\int_{\R^2}\int_{\R^2} \ln(b+|x-y|)F(u_n(x))F(u_n(y))\dd x \dd y\right| \leq C_1\,, \label{bounded:thesis:1}\\
		|A_1(F(u_n),f(u_n)u_n)|=\left|\int_{\R^2}\int_{\R^2} \ln(b+|x-y|)F(u_n(x))f(u_n(y))u_n(y) \dd x \dd y\right| \leq C_2\,.\label{bounded:thesis:2}
	\end{gather}
\end{Lem}
\begin{proof}
	We first show that $\{\nabla u_n\}_n$ is bounded in $L^2(\R^2)$ by combining the information on $\cI$ and its derivative, following the strategy of \cite[Lemma 6.1]{CassaniTarsi2021}. To this aim we introduce the sequence
	\begin{equation*}
		v_n:=\begin{cases}
			\frac{F(u_n)}{f(u_n)}\ \  &\mbox{if}\ \,u_n>0\,,\\
			(1-\tau)u_n & \mbox{if}\ \,u_n<0\,,\\
		\end{cases}
	\end{equation*}
	where $\tau$ is defined in $(f_2)$, and for which, by \eqref{ARcond_delta}, $|v_n|\leq(1-\tau)|u_n|$ hold, hence $\{v_n\}_n\subset L^p_\omega(\R^2)$. Moreover, a simple computation shows that
	\begin{equation*}
		\nabla\left(\frac{F(u_n)}{f(u_n)}\right)=\left(1-\frac{F(u_n)f'(u_n)}{(f(u_n))^2}\right) \nabla u_n\,,
	\end{equation*}
	therefore ($f_2$) implies $|\nabla v_n|\leq C|\nabla u_n|$, from which $\{v_n\}_n\subset D^{1,2}(\R^2)$ and in turn $\{v_n\}_n\subset X$. Therefore, they may be used as test functions for $\cI'(u_n)\in X'$, obtaining
	\begin{equation}\label{ctrl_after_vn}
		\begin{aligned}
			&\int_{\{u_n\geq0\}}|\nabla u_n|^2\left(1-\frac{F(u_n)f'(u_n)}{\left(f(u_n)\right)^2}\right)+(1-\tau)\int_{\{u_n<0\}}|\nabla u_n|^2\\
			&\quad+\int_{\{u_n\geq0\}}\left(\ln|\cdot|\ast F(u_n)\right)F(u_n)+(1-\tau)\int_{\{u_n<0\}}\left(\ln|\cdot|\ast F(u_n)\right)f(u_n)u_n\\
			&=|\cI'(u_n)[v_n]|\leq\|\cI'(u_n)\|_{X'}\|v_n\|\les\|\cI'(u_n)\|_{X'}\|u_n\|=\text{o}_n(1)
		\end{aligned}
	\end{equation}
	by \eqref{un_Cerami}. Since $f\equiv0$ on $\R^-$, the last term in the left-hand side is zero. Combining this with $\cI(u_n)\to c_{mp}$, one infers
	\begin{equation*}
		\begin{aligned}
			&\|\nabla u_n\|_2^2 - 2c_{mp} + \text{o}_n(1)
			=\int_{\R^2}\left(\ln\frac1{|\cdot|}\ast F(u_n)\right)F(u_n)\\
			&=\int_{\{u_n\geq0\}}|\nabla u_n|^2\left(1-\frac{F(u_n)f'(u_n)}{\left(f(u_n)\right)^2}\right)+(1-\tau)\int_{\{u_n<0\}}|\nabla u_n|^2+\text{o}_n(1)\\
			&\leq(1-\tau)\|\nabla u_n\|_2^2+\text{o}_n(1)\,,
		\end{aligned}
	\end{equation*}
	from which
	\begin{equation}\label{un_grad_bdd}
		\begin{aligned}
			\|\nabla u_n\|_2^2 \leq \frac{2c_{mp}}\tau+\text{o}_n(1)\,.
		\end{aligned}
	\end{equation}
	This, together with the first two conditions in \eqref{un_Cerami}, yields
	\begin{equation}\label{bounded:thesis:1_pre}
		|I_0(u_n)|\leq C_1 \qquad\mbox{and}\qquad|A_0(F(u_n),f(u_n)u_n)|\leq C_2\,.
	\end{equation}
	Using \eqref{un_grad_bdd} and \eqref{bounded:thesis:1_pre}, the condition $\cJ(u_n)\to0$ directly implies
	\begin{equation}\label{un_p_bdd}
		\begin{aligned}
			\|u_n\|_p^{2p}\leq C\left(\int_{\R^2}F(u_n)\right)^2\leq C\,,
		\end{aligned}
	\end{equation}
	where the first inequality is due to \eqref{estFbelow}. In light of \eqref{un_p_bdd} and \eqref{un_grad_bdd}, $\{I_2(u_n)\}_n$ is bounded thanks to \eqref{I2:finer:estimate}. Recalling the decomposition $I_0=I_1-I_2$, this and \eqref{bounded:thesis:1_pre} imply \eqref{bounded:thesis:1}. The bound \eqref{bounded:thesis:2} follows by similar arguments using $A_0,\,A_1,\,A_2$. The uniform boundedness of $\|u_n\|_{*,p}$ is then a consequence of Lemma \ref{Lemma3.1tilde}(a), applied with $\varphi_n=F(u_n)$, and \eqref{estFbelow}.
\end{proof}

\begin{Rem}\label{Rmk_nonneg_Cerami}
	Thanks to Lemma \ref{boundedness:Cerami}, from now on we can always suppose that Cerami sequences at level $c_{mp}$ verifying \eqref{un_Cerami} are nonnegative. Indeed, $u_n^-:=\min\{u_n,0\}\in X$ and $u_n^-\leq0$, and thus, recalling that $f\equiv0$ on $\R^-$ by assumption, one has
	\begin{equation*}
		\begin{aligned}
			\|\nabla u_n^-\|_2^2&=\|\nabla u_n^-\|_2^2+\int_{\R^2}\left(\ln|\cdot|\ast F(u_n)\right)f(u_n)u_n^-\\
			&=\cI'(u_n)[u_n^-]\leq\|\cI'(u_n)\|_{X'}\|u_n^-\|=\text{o}_n(1)
		\end{aligned}
	\end{equation*}
	since $\|u_n^-\|\leq\|u_n\|\leq C$ by Lemma \ref{boundedness:Cerami}. This implies that $u_n^-\to0$ in $X$ as $n\to+\infty$ and therefore that $\{u_n^+\}_n$ is a Cerami sequence of $\cI$ at level $c_{mp}$, which henceforth we will simply denote by $\{u_n\}_n$.
\end{Rem}

\subsection{Mountain pass level estimate in the critical case}

Under assumption ($f_{sc}$), the boundedness of the Cerami sequences suffices to carry out the main existence argument, see Section \ref{Section_existence}, since uniform estimates of the nonlinear terms may be deduced by \eqref{estFsc} by choosing a suitably small exponent $\alpha$. This of course is not the case when we are dealing with critical exponential nonlinearities, and we need to prove that under ($f_c$), and in particular taking into account ($f_4$), the critical mountain pass level is below a noncompactness threshold. To this aim, let us introduce the usual Moser sequence $\tw_n:B_\rho\to\R^+$ defined by
\begin{equation*}
	\tw_n(x)=\frac1{\sqrt{2\pi}}
	\begin{cases}
		\sqrt{\ln n}\ \  &\mbox{for}\ \,0\leq|x|\leq\frac \rho n\,,\\
		\\
		\displaystyle \frac{\ln (\rho/|x|)}{\sqrt{\ln n}} & \mbox{for}\ \,\frac\rho n<|x|<\rho\,.
	\end{cases}
\end{equation*}
It is easy to see that $\|\nabla\tw_n\|_2=1$ for all $n\in\N$, and that
\begin{equation*}
	\begin{aligned}
		\|\tw_n\|^p_{\ast,p}
		&=(2\pi)^{1-\frac p2}(\ln n)^{\frac p2}\int_0^{\rho/n}\ln(b+r)\,r\dd r\\
		&\quad +\frac{(2\pi)^{1-\frac p2}}{(\ln n)^{\frac p2}}\int_{\rho/n}^{\rho}\ln^p\left(\tfrac\rho r\right)\,\ln(b+r)\,r\dd r\\
		&\leq(2\pi)^{1-\frac p2}(\ln n)^{\frac p2}\frac{\rho^2\ln(b+\rho/n)}{2n^2}+(2\pi)^{1-\frac p2}\frac{\ln (b+\rho)}{(\ln n)^{\frac p2}}\int_{\rho/n}^{\rho}\ln^p\left(\tfrac\rho r\right)\,r\dd r\,.
	\end{aligned}
\end{equation*}  
The last term in the previous expression can be estimated as follows. On the one hand, for $k\in\N$, 
\begin{equation*}
	\int\ln^k\!\left(\tfrac\rho r\right)\,r\dd r=\frac{r^2}2\sum_{j=0}^k\left(\ln(\tfrac\rho r)\right)^{k-j}\frac{k(k-1)\cdots(k-j+1)}{2^j}\,.
\end{equation*}
On the other hand, since $p$ may be an integer or not, a rough estimate reads as follows 	
\begin{align*}
	\int_{\rho/n}^\rho\ln^p({\textstyle{\tfrac\rho r}})\,r\dd r&\leq \int_{\rho/n}^{\rho}\left\{\ln^{[p]}({\textstyle{\frac \rho r}})+\ln^{[p]+1}({\textstyle{\frac \rho r}})\right\}\,r\dd r\\
	&=\frac{\rho^2[p]!}{2^{[p]+1}}\left(1+\frac{[p]+1}{2}\right)+{\text{o}_n}(1)\,,
\end{align*}
so that eventually 
\begin{equation*}
	1\leq\|\tw_n\|^2\leq 1 +\delta_n + {\text{o}_n}\left(\frac 1{\ln n}\right) ,
\end{equation*}
where 
\begin{eqnarray}\label{delta_n2}
	\delta_n:=\rho^\frac4p\frac{(2\pi)^{\frac 2p-1}}{\ln n}\ln^{\frac 2p}(b+\rho)\left[\frac{[p]!}{2^{[p]+1}}\left(1+\frac{[p]+1}{2}\right)\right]^{\frac 2p}.
\end{eqnarray}      
Hence we normalise the Moser sequence $\{\tw_n\}_n$ by defining
\begin{equation}\label{Moser_seq_normalised}
	w_n:=\frac{\tw_n}{\sqrt{1+\delta_n}}\,,\quad\ n\in\N\,.
\end{equation}

\begin{Lem}\label{MP_level}
	Under ($f_1$), ($f_c$), ($f_2$), ($f_4$), one has 
	\begin{equation}\label{MP_est}
		c_{mp}<\frac{2\pi}{\alpha_0},
	\end{equation}
	where $\alpha_0$ is defined in $(f_c)$.
\end{Lem}
\begin{proof} The proof follows the same arguments of \cite[Lemma 5.2]{CassaniTarsi2021}. We claim that there exists $n$ such that
	\begin{equation}\label{claim}
		\max_{t\geq 0}\cI(tw_n)<\frac{2\pi}{\alpha_0}\,.
	\end{equation}
	Let us argue by contradiction and suppose this is not the case, so that for all $n$ let $t_n>0$ be such that
	\begin{equation}\label{bycontr-assump}
		\max_{t\geq 0}\cI(tw_n)=\cI(t_nw_n)\geq\frac{2\pi}{\alpha_0}\,.
	\end{equation}
	Then $t_n$ satisfies $$\frac{\dd}{\dd t}\bigg|_{t=t_n}\cI(tw_n)=0$$ and
	\begin{equation}\label{t_n^2=}
		t^2_n\geq  \int_{\mathbb R^2}\left[\ln \frac{1}{|\cdot|}\ast F(t_nw_n) \right]t_nw_nf(t_nw_n)\dd x\,,
	\end{equation}
	\begin{equation}\label{t_n^2>}
		t^2_n\geq \frac{4\pi}{\alpha_0} + \int_{\mathbb R^2}\left[\ln\frac1{|\cdot|}\ast F(t_nw_n)\right]F(t_nw_n) \dd x\,.
	\end{equation}
	Note that in \eqref{t_n^2=} we have an inequality instead of the equality since in the energy functional it appears only $\|\nabla w_n\|_2^2\leq\|w_n\|^2=1$.
	\noindent From now on let us suppose $\rho\leq 1/2$. This will simplify a few estimates, since for any $(x,y)\in {\text{supp }}w_n\times {\text{supp }}w_n$ we have $|x-y|>1$, and in turn $\ln(1/|x-y|)>0$. Let us now proceed in three steps.
	
	\vskip0.2truecm
	
	\noindent {\bf{Step 1.}} The following holds: $\limsup_{n \to +\infty} t_n^2 \geq 4\pi/\alpha_0$.
	
	\noindent Let us assume by contradiction that $\limsup_{n}t^2_n<4\pi/\alpha_0$: this implies that, up to a subsequence, there exists a positive constant $\delta_0$ such that $t_n^2\leq 4\pi/\alpha_0-\delta_0$ for $n$ large enough. Since $\rho\leq\frac12$, for any $|x|< \rho$, the set $\{y: |x-y|>1,|y|<\rho\}$ is empty. Recalling that the functions $w_n$ are compactly supported in $B_{\rho}$, we have
	\begin{multline*}
		\int_{\mathbb R^2}\left[\ln \frac{1}{|\cdot|}\ast F(t_nw_n)\right]F(t_nw_n)\dd x\\=\int_{B_{\rho}}\int_{|x-y|\leq 1}\ln\frac1{|x-y|} F(t_nw_n(x))F(t_nw_n(y))\dd x\dd y\geq 0\,,
	\end{multline*}
	a contradiction with \eqref{t_n^2>}.\\
	
	\medskip
	
	\noindent {\bf{Step 2.}} The following holds: $\liminf_{n \to +\infty} t_n^2 \leq4\pi/\alpha_0$\,.
	
	\noindent Let us suppose by contradiction that $\liminf_{n \to +\infty}t_n^2>4\pi/\alpha_0$. Hence, up to a subsequence, there exists a constant $\delta_0>0$ such that
	$$t_n^2\geq \frac{4\pi}{\alpha_0}+\delta_0$$
	as $n \to +\infty$. Let us estimate from below the right hand side of \eqref{t_n^2=} (taking into account the possible negative sign of
	the logarithmic function):
	\begin{equation}\label{est1}
		\begin{aligned}
			\int_{\mathbb R^2}&\left[\ln\frac1{|\cdot|}\ast
			F(t_nw_n)\right]t_nw_nf(t_nw_n)\dd x\\
			&=\int_{\{|x|\leq\frac{\rho}n,\,|y|\leq\frac{\rho}n\}}\ln\frac1{|x-y|}F(t_nw_n(x))t_nw_nf(t_nw_n(y))\dd x\dd y\\
			&\quad+\int_{\R^2 \times \R^2 \setminus\{|x|\leq\frac{\rho}n,\,|y|\leq\frac{\rho}n\}}\ln\frac1{|x-y|}F(t_nw_n(x))t_nw_nf(t_nw_n(y))\dd x\dd y\\
			&=:T_1+T_2\,.
		\end{aligned}
	\end{equation}
	Thanks to $(f_4)$ we have for any $\varepsilon>0$ (here we choose $\varepsilon=\beta/2$),
	\begin{equation*}
		sf(s)F(s)\geq\frac{\beta-\varepsilon}{s^2}\,\e^{2\alpha_0 s^2}=\frac{\beta}{2s^2}\,\e^{2\alpha_0 s^2}, \quad \hbox{ for all
		}s\geq s_{\varepsilon}=s_{\beta}\,.
	\end{equation*}
	By the definition of $w_n$ (see \eqref{Moser_seq_normalised}) and since $|x-y|< 2\rho/n<1$, for $n$ large enough we can estimate $T_1$ as follows
	\begin{equation*}
		\begin{aligned}
			T_1&=\int_{B_{\rho/n}}t_nw_nf(t_nw_n)\left(\int_{B_{\rho/n}}\ln\frac{1}{|x-y|}F(t_nw_n)\dd x\right)\dd y\\
			&= \int_{B_{\rho/n}}t_n\frac{\sqrt{\ln n}}{\sqrt{2\pi(1+\delta_n)}}\,f\left(t_n\frac{\sqrt{\ln n}}{\sqrt{ 2\pi(1+\delta_n)}}\right)\\
			&\quad \times \left(\int_{B_{\rho/n}}\ln\frac{1}{|x-y|}\,F\left(t_n\frac{\sqrt{\ln n}}{\sqrt{2\pi(1+\delta_n)}}\right)\dd x\right)\dd y\\
			&\geq2\pi\beta\frac{\e^{\alpha_0 t_n^2[\pi(1+\delta_n)]^{-1}\ln n}}{2\alpha_0t_n^2 [\pi(1+\delta_n)]^{-1}\ln n}\int_{B_{\rho/n}}\int_{B_{\rho/n}}\ln\frac{1}{|x-y|}\dd x\dd y\,.
		\end{aligned}
	\end{equation*}
	Since
	\begin{equation*}
		\int_{B_{\rho/n}}\int_{B_{\rho/n}}\ln\frac1{|x-y|}\dd x\dd y\geq|B_{\rho/n}|^2\ln\frac{n}{2\rho}=\pi^2\left(\frac{\rho}n\right)^4\ln\frac{n}{2\rho}\,,
	\end{equation*}
	we obtain
	\begin{equation}\label{est1.1}
		T_1\geq \pi^3\rho^4\beta\frac{\e^{(\alpha_0t_n^2[\pi(1+\delta_n)]^{-1}-4) \ln n}}{\alpha_0 t_n^2 [\pi(1+\delta_n)]^{-1} \ln n}\ln\frac{n}{2\rho}\geq \frac{\pi^3\rho^4\beta}{\alpha_0t_n^2}\e^{\left(\frac{\alpha_0}\pi\frac{t_n^2}{1+\delta_n}-4\right)\ln n}
	\end{equation}
	for any $n\geq n(\rho,\beta)$. Note that since $\rho \leq 1/2$ we have
	$$T_2\geq 0\,.$$
	Now, combining \eqref{t_n^2=}, \eqref{est1} and \eqref{est1.1} yields
	\begin{equation}\label{est1.2}
		t_n^4\geq\frac{\pi^3\rho^4\beta}{\alpha_0}\e^{\left(\frac{\alpha_0}\pi\frac{t_n^2}{1+\delta_n}-4\right)\ln  n}
	\end{equation}
	which is a contradiction, either if $t_n\to +\infty$ or $t_n$ stays bounded with $t_n^2\geq\frac{4\pi}{\alpha_0}+\delta_0$. The proof of Step 2 is then complete. Observe that, as a consequence of Step 1 and Step 2
	$$
	t_n^2\to\frac{4\pi}{\alpha_0}\quad \hbox{as } n\to +\infty\,.
	$$
	Moreover, as a byproduct of \eqref{est1.2}, we also have
	$$\e^{\left(\frac{\alpha_0}\pi\frac{t_n^2}{1+\delta_n}-4\right)\ln n}\leq C\,,$$
	for some $C>0$, that is
	\begin{equation*}
		\frac{t_n^2}{1+\delta_n}\leq \frac{4\pi}{\alpha_0}+\frac{C}{\ln n}=\frac{4\pi}{\alpha_0}+{\rm{O}}\left(\frac1{\ln n}\right)\,.
	\end{equation*}
	
	\medskip
	
	\noindent {\bf{Step 3.}} We are now in a position of getting a contradiction and determine the quantity $\mathcal V$ which appears in condition $(f_4)$. We have proved that $t_n^2\to 4\pi/\alpha_0$. Moreover, we also know that $t_n^2\geq 4\pi/\alpha_0$ by \eqref{t_n^2>}. By \eqref{est1.2}, recalling the definition \eqref{delta_n2} of $\delta_n$, we have 
	\begin{equation*}
		\begin{aligned}
			\left(\frac{4\pi}{\alpha_0}\right)^2+\text{o}_n(1)\geq t_n^4 & \geq\frac{\pi^3\rho^4\beta}{\alpha_0}\,\e^{4\left(\frac{\alpha_0}{4\pi}\frac{t_n^2}{1+\delta_n}-1\right)\ln n}\geq\frac{\pi^3\rho^4\beta}{\alpha_0}\,\e^{-4\frac{\delta_n}{1+\delta_n}\ln n}\\
			&\geq\frac{\pi^3\rho^4\beta}{\alpha_0}\,\e^{-\rho^\frac4p(2\pi)^{\frac 2p-1}\ln^{\frac 2p}(b+\rho)\left[\frac{[p]!}{2^{[p]+1}}\left(1+\frac{[p]+1}{2}\right)\right]^{\frac 2p}+\text{o}_n(1)}.
		\end{aligned}
	\end{equation*}
	Passing to the limit, we obtain 
	\begin{equation}\label{contr}
		\frac{16\pi^2}{\alpha_0^2}\geq\frac{\pi^3\rho^4\beta}{\alpha_0}\,\e^{-\rho^\frac4p(2\pi)^{\frac 2p-1}\ln^{\frac 2p}(b+\rho)\left[\frac{[p]!}{2^{[p]+1}}\left(1+\frac{[p]+1}{2}\right)\right]^{\frac 2p}}.
	\end{equation}
	Now set in assumption $(f_4)$ 
	\begin{equation}\label{V}
		\mathcal V:=\inf_{|x|\leq 1/2}\frac{16}{\alpha_0\pi}|x|^{-4}\e^{|x|^\frac4p(2\pi)^{\frac 2p-1}\ln^{\frac 2p}(b+|x|)\left[\frac{[p]!}{2^{[p]+1}}\left(1+\frac{[p]+1}{2}\right)\right]^{\frac 2p}},
	\end{equation}
	a quantity which is actually a minimum, since the right-hand function is continuous and unbounded as $|x|\to 0$ . Finally, since $\beta >\mathcal V$, we can fix $\rho\in (0,1/2]$ such that
	$$
	\beta>\frac{16}{\alpha_0\pi}\rho^{-4}\e^{\rho^\frac4p(2\pi)^{\frac 2p-1}\ln^{\frac 2p}(b+\rho)\left[\frac{[p]!}{2^{[p]+1}}\left(1+\frac{[p]+1}{2}\right)\right]^{\frac 2p}}
	$$
	to get 
	$$
	\frac{\pi^3\rho^4\beta}{\alpha_0}\e^{-\rho^\frac4p(2\pi)^{\frac 2p-1}\ln^{\frac 2p}(b+\rho)\left[\frac{[p]!}{2^{[p]+1}}\left(1+\frac{[p]+1}{2}\right)\right]^{\frac 2p}}>\frac{4\pi}{\alpha_0}\,,
	$$
	which contradicts \eqref{contr} and, therefore, \eqref{bycontr-assump}. This shows that \eqref{claim} holds, and in turn \eqref{MP_est}.
\end{proof}

To avoid trivial solutions, in showing existence we will need to prove a result \textit{à la} Lions (see Section \ref{Section_existence} below). To this end, in the spirit of \cite[Lemma 6.3]{CassaniTarsi2021}, we need to improve the integrability for $F(u_n)$, where $\{u_n\}_n$ is the Cerami sequence given by Lemma \ref{existence:Cerami:sequence}, since this will enable us to uniformly control the terms appearing from an application of the Hardy-Littlewood-Sobolev inequality. Here the mountain pass level estimate given by Lemma \ref{MP_level} plays a crucial role. Unlike \cite[Lemma 6.3]{CassaniTarsi2021}, we cannot apply Ruf's version of the Trudinger-Moser inequality in $H^1(\R^2)$ because of the lack of the mass term in \eqref{SP:equation}; however, since our space $X\subset L^p(\R^2)$, we will make use of the refinement of Cao's inequality in Theorem \ref{GuoLiu_Thm}.

\begin{Lem}\label{Lemma_vn}
	Assume $(f_1)-(f_4)$ and $(f_c)$. Let $\{u_n\}_n \subset X$ be a nonnegative Cerami sequence for $\cI$ at level $c_{mp}<\frac{2\pi}{\alpha_0}$, which is bounded in $X$. Then, the sequence
	\begin{equation*}
		v_n:=G(u_n)\,,\qquad\mbox{where}\quad G(t):=\int_0^t\sqrt{\frac{F(s)f'(s)}{f(s)^2}}\dd s\,,
	\end{equation*}
	has the following properties:
	\begin{enumerate}
		\item[\textit{i})] $\|v_n\|_p\leq\|u_n\|_p$ and $\|v_n\|_{*,p}\leq\|u_n\|_{*,p}$;
		\item[\textit{ii})] $\{v_n\}_n\subset X$ and 
		\begin{equation}\label{vn_cmp}
			\|\nabla v_n\|_2^2=2c_{mp}+\text{o}_n(1)<\frac{4\pi}{\alpha_0}\,;
		\end{equation}
		\item[\textit{iii})] $v_n>\sqrt\tau\,u_n$ a.e. in $\R^2$;
		\item[\textit{iv})] for all $\varepsilon>0$ there exists $t_\varepsilon>0$ such that
		\begin{equation}\label{stima_epsilon_un_vn}
			u_n\leq t_\varepsilon+\frac{v_n}{1-\varepsilon}\quad\ \mbox{for a.e.}\ \ x\in\R^2.
		\end{equation}
	\end{enumerate}
\end{Lem}
\begin{proof}
	Note that $G$ is well-defined and $C^1$ thanks to ($f_2$). By ($f_2$) one has
	\begin{equation*}
		\|\nabla v_n\|_2^2=\int_{\R^2}|\nabla u_n|^2\left(\frac{F(u_n)f'(u_n)}{f(u_n)^2}\right)\leq C\|\nabla u_n\|_2^2\leq C
	\end{equation*}
	and, by $(f_1)$ and H\"older's inequality,
	\begin{equation*}
		\begin{aligned}
			|G(t)|^p&\leq t^{\tfrac p2}\left(\int_0^t\frac{F(s)f'(s)}{f(s)^2}\dd s\right)^{\frac p2}\leq t^{\tfrac p2}\left(-\frac{F(t)}{f(t)}+\lim_{s\to0^+}\frac{F(s)}{f(s)}+t\right)^{\frac p2}\\
			&\leq t^p\left(-\frac{F(t)}{tf(t)}+1\right)^{\frac p2}\leq t^p.
		\end{aligned}
	\end{equation*}
	This yields at once $\|v_n\|_p\leq\|u_n\|_p$ and $\|v_n\|_{*,p}\leq\|u_n\|_{*,p}$. Note also that $G(t)>\sqrt\tau t$, hence $v_n>\sqrt\tau\,u_n$ a.e. in $\R^2$. Recalling now that $u_n\geq0$, and combining \eqref{ctrl_after_vn} with $\cI(u_n)\to c_{mp}$, we infer
	\begin{equation*}
		\begin{aligned}
			2c_{mp}+\text{o}_n(1)&=\|\nabla u_n\|_2^2+\int_{\R^2}\left(\ln|\cdot|\ast F(u_n)\right)F(u_n)\\
			&=\|\nabla u_n\|_2^2-\int_{\R^2}|\nabla u_n|^2\left(1-\frac{F(u_n)f'(u_n)}{f(u_n)^2}\right)=\|\nabla v_n\|_2^2\,.
		\end{aligned}
	\end{equation*}
	Since $c_{mp}<\frac{2\pi}{\alpha_0}$ by assumption, we deduce \eqref{vn_cmp}. Finally, using ($f_3$) and arguing as in \cite[Lemma 6.3]{CassaniTarsi2021}, one proves also \eqref{stima_epsilon_un_vn}.
\end{proof}    

\begin{Lem}\label{Lemma6.3}   
	Assume $(f_1)-(f_4)$ and $(f_c)$. Let $\{u_n\}_n \subset X$ be a nonnegative Cerami sequence for $\cI$ at level $c_{mp}<\frac{2\pi}{\alpha_0}$, which is bounded in $X$. Then, for any $\gamma \in \left[1,\frac{2\pi}{\alpha_0c_{mp}}\right)$ the following uniform bound holds
	\begin{equation*}
		\sup_{n\in\N}\int_{\R^2} (F(u_n))^\gamma \dd x < + \infty\,.
	\end{equation*}
	If, moreover, $u_n\to0$ in $L^s(\R^2)$ for all $s>p$, then for any $\gamma \in \left(1,\frac{2\pi}{\alpha_0c_{mp}}\right)$ one has
	\begin{equation*}
		\|F(u_n)\|_\gamma \to 0 \quad \text{ and } \quad \|f(u_n)u_n\|_\gamma \to 0\,.
	\end{equation*}
\end{Lem}
\begin{proof}
	Since $\{u_n\}_n$ is bounded in the reflexive space $X$, there exists $u\in X$ such that $u_n\rightharpoonup u$ in $X$, which implies that, up to a subsequence, the convergence is strong in $L^q(\R^2)$, for $q \geq p$, and a.e. in $\R^2$ thanks to Proposition \ref{cpt_emb_X}. From Lemma \ref{Lemma_vn}, together with \eqref{estFsc} with $\alpha>\alpha_0$, we may estimate as follows:
	\begin{equation*}
		\begin{aligned}
			\int_{\R^2}|F(u_n)|^\gamma&\leq C_\varepsilon\int_{\{u_n < t_\varepsilon\}}|u_n|^{p\gamma}\\
			&\quad +C_\varepsilon\int_{\{u_n\geq t_\varepsilon\}}\left(t_\varepsilon+\frac{v_n}{1-\varepsilon}\right)^{\gamma r}\e^{\gamma\alpha(1+\varepsilon)\left(t_\varepsilon+\frac{v_n}{1-\varepsilon}\right)^2}\\
			&\leq C_\varepsilon\int_{\R^2}|u_n|^{p\gamma}+C_\varepsilon\int_{\{u_n\geq t_\varepsilon\}}\e^{\gamma\alpha\left(\frac{1+\varepsilon}{1-\varepsilon}\right)^2\|\nabla v_n\|_2^2\left(\frac{v_n}{\|\nabla v_n\|_2}\right)^2}\,,
		\end{aligned}
	\end{equation*}
	where in the last step we used $s^{\gamma r}\leq C_\varepsilon\,\e^{\varepsilon s^2}$ for any $s\geq t_\varepsilon$ and
	\begin{equation}\label{stima_vn_un}
		\left(t_\varepsilon+\frac{v_n}{1-\varepsilon}\right)^2\leq C_\varepsilon t_\varepsilon^2+(1+\varepsilon)\left(\frac{v_n}{1-\varepsilon}\right)^2,
	\end{equation}
	by the $\varepsilon-$Young inequality. Choosing now $\varepsilon$ small enough such that $t_\varepsilon\geq\tau^{-1/2}$, in the set $\{u_n\geq t_\varepsilon\}$ one has $v_n\geq1$, and therefore, by means of the (strict) fine upper bound in \eqref{vn_cmp}, Theorem \ref{GuoLiu_Thm} with $q=p$, $\gamma\in[1,\frac{2\pi}{\alpha_0 c_{mp}})$, and choosing $\alpha$ close to $\alpha_0$ and $\varepsilon$ small, one gets
	\begin{equation}\label{Lem6.3_ultimostep}
		\begin{aligned}
			\int_{\R^2}&|F(u_n)|^\gamma\leq C_\varepsilon\int_{\R^2}|u_n|^{p\gamma}+C_\varepsilon\int_{\R^2}|v_n|^p\e^{\gamma\alpha\left(\frac{1+\varepsilon}{1-\varepsilon}\right)^2\|\nabla v_n\|_2^2\left(\frac{v_n}{\|\nabla v_n\|_2}\right)^2}\\
			&\leq C_\varepsilon\int_{\R^2}|u_n|^{p\gamma}+C_\varepsilon\|\nabla v_n\|_2^p\int_{\R^2}\left|\frac{v_n}{\|\nabla v_n\|_2}\right|^p\e^{\gamma\alpha\left(\frac{1+\varepsilon}{1-\varepsilon}\right)^2\|\nabla v_n\|_2^2\left(\tfrac{v_n}{\|\nabla v_n\|_2}\right)^2}\\
			&\leq C_\varepsilon\|u_n\|_{p\gamma}^{p\gamma}+C_\varepsilon\|v_n\|_p^p\leq C_\varepsilon\|u_n\|_{p\gamma}^{p\gamma}+C_\varepsilon\|u_n\|_p^p\leq C\,,
		\end{aligned}
	\end{equation}
	since $\{u_n\}_n$ is bounded in $X$ and by Proposition \ref{cpt_emb_X}.
	
	Suppose now in addition that $u_n\to0$ in $L^s(\R^2)$ for all $s>p$.  In the same way as in \eqref{Lem6.3_ultimostep} for some $r>0$ and $\alpha>\alpha_0$ one has
	\begin{equation*}
		\begin{aligned}
			&\int_{\R^2}|F(u_n)|^\gamma 
			\leq C_\varepsilon\int_{\R^2}|u_n|^{p\gamma} +C_\varepsilon\int_{\R^2}|v_n|^{r\gamma}\e^{\gamma\alpha\left(\frac{1+\varepsilon}{1-\varepsilon}\right)^2\|\nabla v_n\|_2^2\left(\frac{v_n}{\|\nabla v_n\|_2}\right)^2}\!\dd x\\
			& \quad \leq 
			C_\varepsilon\|u_n\|^{p\gamma}_{p\gamma} + C_\varepsilon\|v_n\|_{\frac r2\sigma'\gamma}^{\frac r2\gamma}\left(\int_{\R^2}|v_n|^{\frac r2\sigma\gamma}\e^{\sigma\gamma\alpha\left(\frac{1+\varepsilon}{1-\varepsilon}\right)^2\|\nabla v_n\|_2^2\left(\frac{v_n}{\|\nabla v_n\|_2}\right)^2}\dd x\right)^{\frac1\sigma}.
		\end{aligned}
	\end{equation*}
	The first term on the right-hand side goes to $0$ as $n \to +\infty$ since $p\gamma>p$; for the second term we proceed by using H\"older's inequality and Theorem \ref{GuoLiu_Thm}: indeed, choosing $\alpha$ close to $\alpha_0$, $\sigma>1$ close to $1$, and $\varepsilon>0$ small, since $\gamma<\frac{2\pi}{c_{mp}\alpha_0}$ and by \eqref{vn_cmp} the exponent is less than $4\pi$; if we then fix $r>2p$, we obtain
	\begin{equation*}
		\begin{aligned}
			\int_{\R^2}|F(u_n)|^\gamma&\les\|u_n\|^{p\gamma}_{p\gamma}+\|v_n\|^{\frac r2\sigma'\gamma}_{\frac r2\gamma}\|\nabla v_n\|_2^{\frac r2\gamma}\les\|u_n\|^{p\gamma}_{p\gamma}+\|v_n\|^{\frac r2\sigma'\gamma}_{\frac r2\gamma}\to0
		\end{aligned}
	\end{equation*}
	as $n \to +\infty$ since $\frac r2\gamma>p\gamma>p$. In the last step we used again \eqref{vn_cmp}.\\
	The proof of $\|f(u_n)u_n\|_\gamma\to0$ follows the same line using \eqref{estFc}.
\end{proof}

\section{Existence in the critical and subcritical case: Proof of Theorem \ref{existence:nontrivial:solution}}\label{Section_existence}
By Lemma \ref{finite:positive:mp:level} we know that $\cI$ satisfies the mountain pass geometry. This yields the existence of a Cerami sequence $\{u_n\}_n \subset X$ for $\cI$ at level $c_{mp}$ defined in \eqref{mountain:pass:level}, satisfying \eqref{un_Cerami} by Lemma \ref{existence:Cerami:sequence}, and which can be assumed nonnegative by Remark \ref{Rmk_nonneg_Cerami}. By Lemma \ref{boundedness:Cerami} such a sequence is bounded in $X$. Suppose by contradiction that $\{u_n\}_n$ is vanishing, namely
\begin{equation}\label{vanishing}
	\liminf_{n \to +\infty}\sup_{y\in\R^2} \int_{B_2(y)}|u_n|^p\dd x = 0\,.
\end{equation}
Since $\{u_n\}_n$ is bounded, then $u_n\to0$ in $L^s(\R^2)$ as $n\to+\infty$ for every $s\in(p,+\infty)$ by \cite[Lemma I.1]{Lions:compact:II1987}. By exploiting the inequality $\ln\left(1+\tfrac bt\right)\les t^{-q}$ for $q\in(0,1]$ and $t\in\R$, and the Hardy-Littlewood-Sobolev inequality, if ($f_c$) is assumed one may estimate
\begin{equation}\label{I2to0}
	\begin{aligned}
		I_2(u_n)&=\intd\intd\ln\left(1+\frac b{|x-y|}\right)F(u_n(x))F(u_n(y))\dd x\dd y\\
		&\les\intd\intd\frac{F(u_n(x))F(u_n(y))}{|x-y|^{\frac{4(\gamma-1)}\gamma}}\dd x\dd y\les\|F(u_n)\|_\gamma^2\to0
	\end{aligned}
\end{equation}
as $n\to+\infty$, which holds for $\gamma\in\left(1,\min\left\{\frac43,\frac{2\pi}{\alpha_0 c_{mp}}\right\}\right)$ by Lemma \ref{Lemma6.3}. On the other hand, when ($f_{sc}$) holds, by combining \eqref{estFsc} with $r=p/2$ and Theorem \ref{GuoLiu_Thm} with $\alpha$ sufficiently small, we get
\begin{equation}\label{sottocritico_gamma}
	\begin{aligned}
		\int_{\R^2}|F(u_n)|^\gamma&\les\intd|u_n|^{\gamma p}+\intd|u_n|^{\gamma\frac p2}\,\e^{\gamma\alpha|u_n|^2}\\
		&\leq\|u_n\|_{\gamma p}^{\gamma p}+\|u_n\|_{\gamma p}^{\gamma\frac p2}\left(\intd|u_n|^{\gamma p}\,\e^{2\gamma\alpha u_n^2}\right)^\frac12\\
		&\leq\|u_n\|_{\gamma p}^{\gamma p}+\|u_n\|_{\gamma p}^{\gamma\frac p2}\|\nabla u_n\|_2^{(\gamma-1)p}\|u_n\|_p^p\to0\,,
	\end{aligned}
\end{equation}
since the last two terms are bounded in $n$, and therefore again $I_2(u_n)\to0$. Analogously, in both cases ($f_c$)-($f_{sc}$), one may show that
\begin{equation*}
	A_2(F(u_n),f(u_n)u_n)=\intd\intd\ln\left(1+\frac b{|x-y|}\right)F(u_n(x))f(u_n(y))u_n(y)\dd x\dd y\to0
\end{equation*}
as $n\to+\infty$. Hence,
\begin{equation}\label{I_0bdd}
	\begin{aligned}
		2c_{mp}+\text{o}_n(1)&=2\cI(u_n)-\cI'(u_n)[u_n]=I_0(u_n)-A_0(F(u_n),f(u_n)u_n)\\
		&=2I_1(u_n)-A_1(F(u_n),f(u_n)u_n)\\
		&\quad -2I_2(u_n)+A_2(F(u_n),f(u_n)u_n)\\
		&=\intd\left(\ln(b+|\cdot|)\ast F(u_n)\right)\left(F(u_n)-f(u_n)u_n\right)\dd x+\text{o}_n(1)<0
	\end{aligned}
\end{equation}
for large $n$ by \eqref{ARcond_delta}, a contradiction. This implies that the vanishing \eqref{vanishing} does not occur. Consequently, there exist $\delta>0$ and a sequence $\{y_n\}_n\subset\R^2$ such that (up to a subsequence)
\begin{equation*}
	\int_{B_1(y_n)}|u_n|^p \dd x>\delta\,.
\end{equation*}
Defining $\tu_n:=u_n(\cdot+y_n)$, it is easy to see that
\begin{equation}\label{tu_n_nontrivial}
	\int_{B_1}|\tu_n|^p\dd x>\delta\,,
\end{equation}
and $\|\nabla\tu_n\|_2+\|\tu_n\|_p=\|\nabla u_n\|_2+\|u_n\|_p\leq C$. Therefore $\{\tu_n\}_n$ is a bounded sequence in $D^{1,2}(\R^2)\cap L^p(\R^2)$. Moreover,
$$\|\tu_n\|_{\ast,p}^p=\int_{\R^2}\ln(b+|x-y_n|)|u_n(x)|^p\dd x\leq\|u_n\|_{\ast,p}^p+\ln(b+|y_n|)\|u_n\|_p^p$$
by \eqref{log_sum}, thus $\tu_n\in X$ for all $n\in\N$. Since $I_i(\tu_n)=I_i(u_n)$ for all $i\in\{0,1,2\}$, we also deduce that $\{\tu_n\}_n$ satisfies
\begin{equation}\label{Cerami_tu_n}
	\cI(\tu_n)\to c_{mp}, \quad \cI'(\tu_n)[\tu_n] \to 0, \quad \cJ(\tu_n)\to0\,.
\end{equation}
Note that second condition in \eqref{Cerami_tu_n} is weaker than the corresponding in Lemma \ref{boundedness:Cerami}. Furthermore, since $\left|I_0(\tu_n)\right|\leq C$ as in \eqref{I_0bdd} and $I_2(\tu_n)\leq C$ as in \eqref{I2to0}, we conclude that $I_1(\tu_n)\leq C$. Applying Lemma \ref{Lemma3.1tilde}(a) with $\varphi_n=F(\tu_n)$ we get
\begin{equation}\label{bound_unif_norm_log}
	\|F(\tu_n)\|_{*,1}\leq C\qquad\mbox{and thus}\qquad\|\tu_n\|_{\ast,p}\leq C
\end{equation}
by \eqref{estFbelow}. Since $\{\tu_n\}_n$ is hence bounded in $X$, there exists $\tu\in X$ such that $\tu_n\rightharpoonup\tu$ in $X$ and $\tu_n\to\tu$ in $L^q(\R^2)$ for all $q\geq p$ as well as a.e. in $\R^2$ by Proposition \ref{cpt_emb_X}. By \eqref{tu_n_nontrivial} it is easy to see that
\begin{equation}\label{u_nontrivial}
	\delta<\|\tu_n\|_{L^p(B_1)}^p\to\|\tu\|_{L^p(B_1)}^p\,,
\end{equation}
hence $\tu\not\equiv0$.  We next show that $\cI'(\tu_n)\to0$ in $X'$ which, together with \eqref{Cerami_tu_n}, makes $\tu_n$ verify the same properties as $u_n$ in Lemma \ref{boundedness:Cerami}. Indeed, $\cI(\tu_n)[\varphi]=\cI(u_n)[\varphi(\cdot-y_n)]$ for all $\varphi\in X$ and
\begin{equation}\label{varphi-yn}
	\|\varphi(\cdot-y_n)\|^2\leq\|\varphi\|^2+\ln(b+|y_n|)\|\varphi\|_p^p\,.
\end{equation}
Moreover, if $|y_n|\leq b$, then $\ln(b+|y_n|)\leq\ln(2b)$, while if $|y_n|>b$ one can use \eqref{tu_n_nontrivial} to show that
\begin{equation}\label{estimate_yn}
	\begin{aligned}
		\|u_n\|_{\ast,p}^p&\geq\int_{B_1}\ln(b+|x+y_n|)|\tu_n(x)|^p\dd x\geq\ln|y_n|\int_{B_1}|\tu_n|^p\\
		&\geq\delta\ln|y_n|\geq\delta\ln(b+|y_n|)\frac{\ln b}{\ln(2b)}\,.
	\end{aligned}
\end{equation}
In the third step we used the simple inequality $|x+y_n|\geq|y_n|-|x_n|\geq|y_n|-1\geq|y_n|-b$, while the last step follows from $\frac{\ln|y_n|}{\ln(b+|y_n|)}>\frac{\ln b}{\ln(2b)}$, which holds by monotonicity. Eventually from \eqref{varphi-yn}-\eqref{estimate_yn} we get
\begin{equation*}
	\|\varphi(\cdot-y_n)\| \leq \left[\|\varphi\|^2+\left(\frac1\delta\frac{\ln(2b)}{\ln b}\|u_n\|_{\ast,p}^p+\ln(2b)\right)\|\varphi\|_p^p\right]^\frac12,
\end{equation*}
which implies that
\begin{equation*}
	\begin{aligned}
		|\cI'(\tu_n)[\varphi]|&\leq\|\cI'(u_n)\|_{X'}\left[\|\varphi\|^2+\ln(2b)\left(\frac1{\delta\ln b}\|u_n\|_{\ast,p}^p+1\right)\|\varphi\|_p^p\right]^\frac12\\
		&\leq C\|\cI'(u_n)\|_{X'}\left[\|\varphi\|^2+\|\varphi\|_p^p\right]\to0\,,
	\end{aligned}
\end{equation*}
as $n\to+\infty$, by recalling that $\|\cI'(u_n)\|_{X'}\to0$ and that $\|u_n\|_{\ast,p}$ is uniformly bounded  by Lemma \ref{boundedness:Cerami}.\\
Since then $\cI'(\tu_n)\to0$ in $X'$ and $\cI'(\tu_n)[\tu_n] \to 0$, we deduce that
\begin{equation}\label{final_equality}
	\begin{aligned}
		\text{o}_n(1)&=\cI'(\tu_n)[\tu_n-\tu]\\
		&=\int_{\R^2}\nabla\tu_n \cdot \nabla(\tu_n-\tu)+A_0(F(\tu_n),f(\tu_n)(\tu_n-\tu))\,.
	\end{aligned}
\end{equation}
First,
\begin{equation}\label{final_equality_gradients}
	\begin{aligned}
		\int_{\R^2}\nabla\tu_n \cdot \nabla(\tu_n-\tu)&=\int_{\R^2}|\nabla(\tu_n-\tu)|^2+\int_{\R^2}\nabla\tu\cdot\nabla(\tu_n-\tu)\\
		&=\|\tu_n-\tu\|^2_{D^{1,2}(\R^2)}+\text{o}_n(1)
	\end{aligned}
\end{equation} 
since $\tu_n\rightharpoonup\tu$ in $D^{1,2}(\R^2)$. Moreover, we claim that
\begin{equation}\label{final_equality_log_good}
	A_2(F(\tu_n),f(\tu_n)(\tu_n-\tu))\les\|F(\tu_n)\|_\gamma\|f(\tu_n)(\tu_n-\tu)\|_\gamma\to0\,.
\end{equation}
Indeed, if ($f_{sc}$) is assumed, it is easy to deal with both terms and prove that the first is uniformly bounded and that the second converges to $0$, with analogous estimates as done e.g. in \eqref{sottocritico_gamma}. On the other hand, if ($f_c$) holds, $\|F(\tu_n)\|_\gamma\leq C$ by Lemma \ref{Lemma6.3}, while we can estimate the second term as follows, by defining $v:=G(\tu_n)$ and by means of Lemma \ref{Lemma_vn}:
\begin{equation*}
	\begin{aligned}
		\intd&\left|f(\tu_n)(\tu_n-\tu)\right|^\gamma\leq C\|\tu_n\|_{p\gamma}^{(p-1)\gamma}\|\tu_n-\tu\|_{p\gamma}^\gamma+C_\varepsilon\|\tu_n-\tu\|_{\sigma'\gamma}^\gamma\\
		&\quad\times\|\nabla v_n\|_2^{r\gamma}\left(\intd\left(\frac{|v_n|}{\|\nabla v_n\|_2}\right)^{r\gamma\sigma}\e^{\gamma\alpha\sigma\left(\frac{1+\varepsilon}{1-\varepsilon}\right)^2\|\nabla v_n\|_2^2\left(\frac{v_n}{\|\nabla v_n\|_2}\right)^2}\right)^\frac1\sigma.
	\end{aligned}
\end{equation*}
The first term tends to $0$ as $n\to+\infty$, since $\{\tu_n\}_n$ is bounded in $X$ and by the continuous and compact embeddings in Lebesgue spaces provided by Proposition \ref{cpt_emb_X}, since $\gamma>1$. By choosing $\gamma>1$ and $\sigma>1$ both close to $1$, $\alpha>\alpha_0$ close to $\alpha_0$, $\varepsilon>0$ small, and $r=p$, by \eqref{vn_cmp} and Theorem \ref{GuoLiu_Thm}, this yields
\begin{equation*}
	\begin{aligned}
		\intd&\left|f(\tu_n)(\tu_n-\tu)\right|^\gamma\les \text{o}_n(1)+\|\tu_n-\tu\|_{\sigma'\gamma}^\gamma\|\nabla v_n\|_2^{p\left(\gamma-\frac1\sigma\right)}\|v_n\|_p^\frac p\sigma\\
		&\les \text{o}_n(1)+\left(\frac{4\pi}{\alpha_0}\right)^{\frac p2\left(\gamma-\frac1\sigma\right)}\|\tu_n-\tu\|_{\sigma'\gamma}^\gamma\|\tu_n\|_p^{\frac p\sigma}=\text{o}_n(1)\,,
	\end{aligned}
\end{equation*}
since $\gamma-\frac1\sigma>0$ and again by Proposition \ref{cpt_emb_X}.
Therefore, combining \eqref{final_equality}, \eqref{final_equality_gradients}, and \eqref{final_equality_log_good} we infer
\begin{equation}\label{final_equality_post}
	\begin{aligned}
		\text{o}_n(1)
		&=\|\tu_n-\tu\|^2_{D^{1,2}}+A_1(F(\tu_n),f(\tu_n)(\tu_n-\tu))\\
		&=\|\tu_n-\tu\|^2_{D^{1,2}}+A_1(F(\tu_n),\tf(\tu_n-\tu)(\tu_n-\tu))\\
		&\quad+A_1\left(F(\tu_n),\left(\tf(\tu_n)-\tf(\tu_n-\tu)\right)(\tu_n-\tu)\right)
	\end{aligned}   
\end{equation}
as $n\to+\infty$, where we have defined
\begin{equation}\label{tf}
	\tf(t)=\left\{
	\begin{array}{ll}
		f(t) & \mbox{for }t\geq0\,,\\
		-f(-t) & \mbox{for }t<0\,,
	\end{array}\right.
\end{equation}
since $u_n\geq0$. Note that if we prove that
\begin{equation}\label{convergence_aim}
	A_1(F(\tu_n),\tf(\tu_n-\tu)(\tu_n-\tu))\to0\,,
\end{equation}
one could apply Lemma \ref{Lemma3.1tilde}(b) and obtain $\|\tf(\tu_n-\tu)(\tu_n-\tu)\|_{\ast,1}\to0$, which in turn implies that $\|\tu_n-\tu\|_{\ast,p}\to0$ as $n\to+\infty$ by ($f_{sc}$) or ($f_c$). Together with $\tu_n\to\tu$ in $D^{1,2}(\R^2)$, this would imply that $\tu_n\to\tu$ in $X$, which concludes the proof since $\{\tu_n\}_n$ is a Cerami sequence and $\cI$ is of class $C^1$. Since the first two terms in the right-hand side of \eqref{final_equality_post} are positive, we are then lead to show that
\begin{equation}\label{convergence_aim_2}
	A_1\left(F(\tu_n),\left(\tf(\tu_n)-\tf(\tu_n-\tu)\right)(\tu_n-\tu)\right)\to0
\end{equation}
as $n\to+\infty$. Using \eqref{log_sum}, we obtain
\begin{equation}\label{convergence_aim_2_step1}
	\begin{aligned}
		\bigg|\int_{B_R}&\!\big(\ln(b+|\cdot|)\ast F(\tu_n)\big)\left[\left(\tf(\tu_n)-\tf(\tu_n-\tu)\right)(\tu_n-\tu)\right]\bigg|\\
		&\leq\|F(\tu_n)\|_{\ast,1}\intd\left|\tf(\tu_n)-\tf(\tu_n-\tu)\right||\tu_n-\tu| \\
		&\quad+\|F(\tu_n)\|_1\intd\ln(b+|x|)\left|\tf(\tu_n)-\tf(\tu_n-\tu)\right||\tu_n-\tu|\,.
	\end{aligned}
\end{equation}
First, note that $\|F(\tu_n)\|_1\les\|F(\tu_n)\|_{\ast,1}\leq C$ by \eqref{bound_unif_norm_log}. If ($f_{sc}$) holds, then by \eqref{estFsc} with $r=p$, one infers
\begin{align}\label{f(u_n)_sc}
	\nonumber		\intd f(\tu_n)|\tu_n-\tu|&\les\intd|\tu_n|^{p-1}|\tu_n-\tu|+\intd|\tu_n|^{p-1}\e^{\alpha\tu_n^2}\,|\tu_n-\tu|\\
	&\leq\|\tu_n\|_p^{p-1}\|\tu_n-\tu\|_p+\left(\intd|\tu_n|^p\,\e^{\alpha\frac p{p-1}\tu_n^2}\right)^\frac{p-1}p\!\|\tu_n-\tu\|_p\\
	\nonumber	&\les\|\tu_n\|_p^{p-1}\|\tu_n-\tu\|_p=\text{o}_n(1)\,,
\end{align}
by Theorem \ref{GuoLiu_Thm} taking $\alpha<4\pi\frac{p-1}p\left(\sup_n\|\nabla\tu_n\|_2\right)^{-2}$. By a similar estimate, one may also show that $\intd\tf(\tu_n-\tu)|\tu_n-\tu|=\text{o}_n(1)$. 

\noindent On the other hand, if $f$ is of critical growth, namely ($f_c$) holds, defining $\tv_n=G(\tu_n)$, by Lemma \ref{Lemma_vn} one finds
\begin{equation}\label{f(u_n)}
	\begin{aligned}
		&\intd f(\tu_n)|\tu_n-\tu|\\
		&\les\intd|\tu_n|^{p-1}|\tu_n-\tu| + \intd\left(t_\varepsilon+\tfrac{\tv_n}{1-\varepsilon}\right)^{r-1}\e^{\alpha\left(t_\varepsilon+\frac{\tv_n}{1-\varepsilon}\right)^2}|\tu_n-\tu|\\
		&\les\|\tu_n\|_p^{p-1}\|\tu_n-\tu\|_p +\!\left(\!\intd\!\left(t_\varepsilon+\tfrac{\tv_n}{1-\varepsilon}\right)^{(r-1)\sigma'}\!\!\e^{\alpha\sigma'\left(t_\varepsilon+\frac{\tv_n}{1-\varepsilon}\right)^2}\right)^\frac1{\sigma'}\!\|\tu_n-\tu\|_\sigma\\
		&=\text{o}_n(1)\,,
	\end{aligned}
\end{equation}
where $\sigma'>1$ is close to $1$ and $\alpha>\alpha_0$ close to $\alpha_0$, so that the last term is bounded, and by using that $\tu_n\to\tu$ in $L^q(\R^2)$ for $q\geq p$. Moreover, by \eqref{estFc} and the $\varepsilon-$Young inequality,
\begin{align}
	\label{f(u_n-u)}
	&\intd\big|\tf(\tu_n-\tu)\big||\tu_n-\tu|\les\intd|\tu_n-\tu|^p+\intd|\tu_n-\tu|^r\e^{\alpha|\tu_n-\tu|^2} \nonumber\\
	&\leq \text{o}_n(1)+\left(\intd|\tu_n-\tu|^r\e^{(1+\varepsilon)\alpha s\,\tu_n^2}\right)^\frac1s\left(\intd|\tu_n-\tu|^r\e^{C_\varepsilon\alpha s'\,\tu^2}\right)^\frac1{s'} \nonumber\\
	&\leq\text{o}_n(1)+\left(C\!\int_{\{\tu_n\leq1\}}\!|\tu_n-\tu|^r+\int_{\{\tu_n>1\}}\!|\tu_n-\tu|^r|\tu_n|^\kappa\e^{(1+\varepsilon)\alpha s\,\tu_n^2}\right)^\frac1s\\
	&\quad\times\left(C\int_{\{\tu\leq1\}}|\tu_n-\tu|^r+\int_{\{\tu>1\}}|\tu_n-\tu|^r|\tu|^\kappa\e^{C_\varepsilon\alpha s'\,\tu^2}\right)^\frac1{s'} \nonumber\\
	&\leq \text{o}_n(1)+\left(\text{o}_n(1)+\|\tu_n-\tu\|_{rq'}^{r/s}\left(\intd|\tu_n|^{\kappa q}\e^{(1+\varepsilon)\alpha sq\,\tu_n^2}\right)^\frac1{sq}\right) \nonumber\\
	&\quad\times\left(\text{o}_n(1)+\|\tu_n-\tu\|_{rq'}^{r/s'}\left(\intd|\tu|^{\kappa q}\e^{C_\varepsilon\alpha s'\,{\tu^2}}\right)^\frac1{s'q}\right). \nonumber
\end{align}
Defining $v_n:=G(\tu_n)$, by Lemma \ref{Lemma_vn} we can estimate the first remaining term as follows:
\begin{equation*}
	\begin{aligned}
		\intd|\tu_n|^{\kappa q}\e^{(1+\varepsilon)\alpha sq\,\tu_n^2}&\leq\frac1{\tau^{\kappa/2}}\intd|v_n|^{\kappa q}\e^{(1+\varepsilon)\alpha sq\left(t_\varepsilon+\frac{v_n}{1-\varepsilon}\right)^2}\\
		&\leq\frac1{\tau^{\kappa/2}}\e^{(1+\varepsilon)\alpha sC_\varepsilon t_\varepsilon^2}\intd|v_n|^{\kappa q}\e^{\left(\frac{1+\varepsilon}{1-\varepsilon}\right)^2\alpha sqv_n^2}.
	\end{aligned}
\end{equation*}
Fixing $\varepsilon>0$ close to $0$, $\alpha>\alpha_0$ close to $\alpha_0$, $s>1$ close to $1$, and $\kappa=p/q$, recalling \eqref{vn_cmp}, we can now use Theorem \ref{GuoLiu_Thm} and find that
\begin{equation*}
	\begin{aligned}
		\intd|\tu_n|^{\kappa q}\e^{(1+\varepsilon)\alpha sq\,\tu_n^2}&\leq C\|v_n\|_p^p\leq C\|u_n\|_p^p\leq C
	\end{aligned}
\end{equation*}
since $\{\tu_n\}_n$ is bounded in $L^p(\R^2)$. On the other hand, the second integral on the right in \eqref{f(u_n-u)} is independent on $n$, and therefore, for the above choices of the parameters it remains bounded again by Theorem \ref{GuoLiu_Thm}. Consequently, combining \eqref{f(u_n-u)} with \eqref{f(u_n)_sc} or \eqref{f(u_n)}, we get 
\begin{equation}\label{f(u_n-u)-f(u_n)}
	\intd\big|\tf(\tu_n)-\tf(\tu_n-\tu)\big||\tu_n-\tu|=\text{o}_n(1)\,.
\end{equation}
It remains to prove that also the last term in \eqref{convergence_aim_2_step1} vanishes. To this aim, since $f\in C^1(\R)$ and so is $\tf$ (cf. \eqref{tf}), by Lagrange's theorem,
\begin{equation*}
	\tf(\tu_n)-\tf(\tu_n-\tu)=\tf'(w_n)\tu=f'(|w_n|)\tu\,,
\end{equation*}
where $w_n:=\theta_n\tu_n+(1-\theta_n)\tu$ with $\theta_n:\R^2\to[0,1]$. Therefore, splitting $\R^2=B_R\cup B_R^c$ for a fixed $R>0$, we obtain
\begin{align}
	\label{convergence_aim_2_step2}
	\intd\ln(b+|x|)&\left|\tf(\tu_n)-\tf(\tu_n-\tu)\right||\tu_n-\tu|\dd x \nonumber\\
	&\ \ \leq\ln(b+R)\!\int_{B_R}\!\big|\tf(\tu_n)-\tf(\tu_n-\tu)\big||\tu_n-\tu| \dd x\\
	&\quad +\!\int_{B_R^c}\!\!\ln(b+|x|)|f'(|w_n|)|\,\tu\,|\tu_n-\tu|\dd x \nonumber\\
	& \leq \text{o}_n(1)+\int_{B_R^c}\!\!\ln(b+|x|)|f'(|w_n|)|\,\tu\,|\tu_n-\tu|\dd x \nonumber
\end{align}
by \eqref{f(u_n-u)-f(u_n)}. To estimate this last term, we make use of assumption ($f_5$) on $f'$, which implies that
$$f'(t)\leq C_1|t|^{p-2}+C_2|t|^r\e^{\alpha t^2}$$
for $r>0$ and $\alpha>\alpha_0$, and some constants $C_1,\,C_2>0$. Therefore we may estimate as follows
\begin{equation}\label{T_1+T_2}
	\begin{aligned}
		\int_{B_R^c}\!\!\ln(b&+|x|)|f'(|w_n|)|\,\tu\,|\tu_n-\tu|\dd x\\
		&\les \int_{B_R^c}\!\!\ln(b+|x|)|w_n|^{p-2}\,\tu|\tu_n-\tu|\dd x\\
		&\quad+\int_{B_R^c}\!\!\ln(b+|x|)|w_n|^r\e^{\alpha w_n^2}\,\tu|\tu_n-\tu|=:T_1+T_2\,,
	\end{aligned}
\end{equation}
and consider separately the two integrals. By H\"older's inequality, using $|w_n|\leq|\tu_n|+|\tu|$, and $\|\tu_n\|_{\ast,p}\leq C$ independently of $n$, we first have
\begin{equation}\label{T1}
	\begin{aligned}
		T_1&\les\!\left(\int_{B_R^c}\!\!\ln(b+|x|)|\tu|^p\dd x\right)^\frac1p\!\|\tu_n-\tu\|_{\ast,p}\!\left(\intd\!\ln(b+|x|)(|\tu_n|^p+|\tu|^p)\dd x\right)^\frac{p-2}p\\
		&\les\!\left(\int_{B_R^c}\!\!\ln(b+|x|)|\tu|^p\dd x\right)^\frac1p\left(\|\tu_n\|_{\ast,p}+\|\tu\|_{\ast,p}\right)^{p-1}\\
		&\leq C\left(\int_{B_R^c}\!\!\ln(b+|x|)|\tu|^p\dd x\right)^\frac1p\to0
	\end{aligned}
\end{equation}
as $R\to+\infty$ by Lebesgue's dominated convergence theorem, since $\|\tu\|_{*,p}<+\infty$. Next, again by H\"older's inequality,
\begin{equation}\label{T2}
	\begin{aligned}
		T_2&\leq\left(\intd\ln(b+|x|)|w_n|^{r\sigma'}\e^{\alpha\sigma'w_n^2}\dd x\right)^\frac1{\sigma'}\!
		\left(\int_{B_R^c}\!\!\ln(b+|x|)|\tu|^{s\sigma}\dd x\right)^\frac1{s\sigma}\\
		&\quad\times\left(\intd\ln(b+|x|)|\tu_n-\tu|^{s'\sigma}\dd x\right)^\frac1{s'\sigma}.
	\end{aligned}
\end{equation}
We first note that 
\begin{equation}\label{T2_III}
	\left(\intd\ln(b+|x|)|\tu_n-\tu|^{s'\sigma}\dd x\right)^\frac1{s'\sigma}\les\|\tu_n\|+\|\tu\|\leq C
\end{equation}
by Proposition \ref{cpt_emb_X_log}. Moreover, as in \eqref{T1},
\begin{equation}\label{T2_II}
	\int_{B_R^c}\!\!\ln(b+|x|)|\tu|^{s\sigma}\dd x\to0
\end{equation}
as $R\to+\infty$, again by Proposition \ref{cpt_emb_X_log}. Let us now focus on the first term in \eqref{T2}, which we denote by $T_2^1$, aiming at an estimate independent of $n$. An $\varepsilon-$Young's inequality in the exponent yields
\begin{equation}\label{T2_I}
	\begin{aligned}
		\big(T_2^1\big)^{\sigma'}&\les\intd\ln(b+|x|)\left(|\tu_n|^{r\sigma'}+|\tu|^{r\sigma'}\right)\e^{\alpha\sigma'\left[(1+\varepsilon)\tu_n^2+C_\varepsilon\tu^2\right]}\dd x\\
		&\leq\left(\intd\!\ln(b+|x|)|\tu_n|^{r\sigma'}\e^{\alpha\sigma'\nu(1+\varepsilon)\tu_n^2}\dd x\right)^\frac1\nu\!\!\left(\intd\!\ln(b+|x|)|\tu_n|^{r\sigma'}\e^{\alpha\sigma'\nu'C_\varepsilon\tu^2}\dd x\right)^\frac1{\nu'}\\
		&\quad+\left(\intd\!\ln(b+|x|)|\tu|^{r\sigma'}\e^{\alpha\sigma'\overline\nu(1+\varepsilon)\tu_n^2}\dd x\right)^\frac1{\overline\nu} \!\!\left(\intd\!\ln(b+|x|)|\tu|^{r\sigma'}\e^{\alpha\sigma'\overline\nu'C_\varepsilon\tu^2}\dd x\right)^\frac1{\overline\nu'}\\
		&=:S_1S_2+S_3S_4\,.
	\end{aligned}
\end{equation}
By Theorem \ref{thm_wTp}, we first have
\begin{equation}\label{S4}
	S_4<+\infty\,.
\end{equation}
Moreover, defining again $v_n:=G(\tu_n)$, by Lemma \ref{Lemma_vn} and \eqref{stima_vn_un} we estimate $S_1$ as follows:
\begin{equation*}
	\begin{aligned}
		S_1^\nu&\leq\frac{\e^{\alpha\sigma'\nu(1+\varepsilon)C_\varepsilon t_\varepsilon^2}}{\tau^{r\sigma'/2}}\intd\ln(b+|x|)|v_n|^{r\sigma'}\e^{\alpha\sigma'\nu\frac{(1+\varepsilon)^2}{(1-\varepsilon)^2}v_n^2}\dd x\\
		&=C_\varepsilon\|\nabla v_n\|_2^{r\sigma'}\!\intd\ln(b+|x|)\left(\frac{|v_n|}{\|\nabla v_n\|_2}\right)^{r\sigma'}\!\e^{\alpha\sigma'\nu\left(\frac{1+\varepsilon}{1-\varepsilon}\right)^2\|\nabla v_n\|_2^2\left(\frac{v_n}{\|\nabla v_n\|_2}\right)^2}\dd x\,.
	\end{aligned}
\end{equation*}
Choose $\alpha>\alpha_0$ close to $\alpha_0$, $\sigma',\,\nu>1$ both close to $1$, and $\varepsilon$ small, so that
$$\alpha\sigma'\nu\left(\frac{1+\varepsilon}{1-\varepsilon}\right)^2\|\nabla v_n\|_2^2<4\pi$$
by \eqref{vn_cmp}, and then $r=p/\sigma'$. Since $\|v_n\|_{\ast,p}\leq\|\tu_n\|_{\ast,p}\leq C$, we can apply the weighted Cao's inequality of Theorem \ref{propCaoweighted} and get
\begin{equation}\label{S1}
	S_1^\nu\leq C(\|v_n\|_{\ast,p})\,\|\nabla v_n\|_2^p\leq C
\end{equation}
again by \eqref{vn_cmp}. It remains only to bound the terms $S_2$ and $S_3$, which mix $\tu_n$ and $\tu$. On the one hand,
\begin{equation}\label{S3}
	\begin{aligned}
		S_3^{\overline\nu}&\leq\int_{\{\tu\leq\tu_n\}}\ln(b+|x|)|\tu_n|^p\e^{\alpha\sigma'\overline\nu(1+\varepsilon)\tu_n^2}\dd x\\
		&\quad+\int_{\{\tu>\tu_n\}}\ln(b+|x|)|\tu|^p\e^{\alpha\sigma'\overline\nu(1+\varepsilon)\tu^2}\dd x\leq C
	\end{aligned}
\end{equation}
by choosing $\overline\nu>1$ close to $1$ and reasoning as for \eqref{S1} for the first term, and by Theorem \ref{thm_wTp} for the second. On the other hand,
\begin{equation}\label{S2}
	\begin{aligned}
		S_2^{\nu'}&=\intd\ln(b+|x|)|\tu_n|^p\e^{\alpha\sigma'\nu'C_\varepsilon\tu^2}\dd x\\
		&\leq C_\varepsilon\int_{\{\tu\leq1\}}\ln(b+|x|)|\tu_n|^p\dd x+\int_{\{\tu>1\}}\ln(b+|x|)|\tu_n|^p\tu\,\e^{\alpha\sigma'\nu'C_\varepsilon\tu^2}\dd x\\
		&\leq C+\left(\intd\ln(b+|x|)|\tu_n|^{pp'}\dd x\right)^\frac1{p'}\!\left(\intd\ln(b+|x|)\,\tu^p\,\e^{\alpha\sigma'\nu'pC_\varepsilon\tu^2}\dd x\right)^\frac1p\\
		&\leq C\,,
	\end{aligned}
\end{equation}
by Proposition \ref{cpt_emb_X_log} and Theorem \ref{thm_wTp}. Hence, combining \eqref{S4}-\eqref{S2}, we obtain that the term $T_2^1$ in \eqref{T2_I} is uniformly bounded. This, together with \eqref{T2_III} and \eqref{T2_II}, imply that $T_2\to0$ as $R\to+\infty$. Recalling that $T_1$ had the same behaviour by \eqref{T1}, by \eqref{T_1+T_2} and \eqref{convergence_aim_2_step2} this yields
$$\intd\ln(b+|x|)\left|\tf(\tu_n)-\tf(\tu_n-\tu)\right||\tu_n-\tu|\dd x\to0$$
as $n\to+\infty$, which leads to \eqref{convergence_aim_2}. By \eqref{final_equality_post}, this eventually implies both $\|\nabla\tu_n-\nabla\tu\|_2\to0$ and \eqref{convergence_aim}. By Lemma \ref{Lemma3.1tilde}(b) then $\|\tu_n-\tu\|_{\ast,p}\to0$ as $n\to+\infty$ and therefore $\tu_n\to\tu$ in $X$. Since $\cI$ is a $C^1$-functional, then $\tu$ is a weak solution of \eqref{SP:system}, which is nontrivial thanks to \eqref{u_nontrivial}. Since $\tu_n\geq0$, by Remark \ref{Rem_pos_sol} the solution $\tu$ is positive in $\R^2$.

\section{Back to the system: Proof of Theorem \ref{Thm_SP}}\label{Section_SPsystem}

Let $u\in X$ be the weak solution of the Choquard equation \eqref{SP:equation} given by Theorem \ref{existence:nontrivial:solution} and define
$$\Phi_u(x):=\intd\ln\left(\frac1{|x-y|}\right)F(u(y))\dd y\,.$$
Following the approach of \cite{BCT}, we aim at proving that $\Phi_u$ is a solution of the system \eqref{SP:system} in the sense of Definition \ref{sol_SP}. First, we show that $\Phi_u\in L_s(\R^2)$, for all $s>0$:
\begin{equation*}
	\begin{aligned}
		\intd&\frac{|\Phi_u(x)|}{1+|x|^{2+2s}}\dd x\leq\intd F(u(y))\left(\intd\left|\ln\frac1{|x-y|}\right|\frac1{1+|x|^{2+2s}}\dd x\right)\dd y\\
		&\leq\intd F(u(y))\left(\int_{\{|x-y|>1\}}\frac{\ln|x-y|}{1+|x|^{2+2s}}\dd x+\int_{\{|x-y|\leq1\}}\ln\frac1{|x-y|}\dd x\right)\dd y\\
		&\leq\|F(u(y))\|_1\intd\frac{\ln(1+|x|)}{1+|x|^{2+2s}}\dd x+\|F(u(y))\|_{\ast,1}\intd\frac{\dd x}{1+|x|^{2+2s}}\\
		&\quad+\|\ln(\cdot)\|_{L^1(B_1)}\intd F(u(y))\dd y<+\infty
	\end{aligned}
\end{equation*}
for all $s>0$, using \eqref{log_sum} with $b=1$, and Theorem \ref{thm_wTp} since $u\in X$. Define now the function
$$\tw_u(x):=\intd\ln\left(\frac{1+|y|}{|x-y|}\right)F(u(y))\dd y\,,$$
which we know by \cite[Lemma 2.3]{H} to be a solution in the sense of Definition \ref{sol_Poisson} of $-\Delta\tw_u=\ff$ in $\R^2$, where $\ff:=F(u)\in L^1(\R^2)$, and compute the difference 
\begin{equation*}
	\begin{aligned}
		\tw_u(x)-\Phi_u(x)&=\intd\left(\ln\left(\frac{1+|y|}{|x-y|}\right)-\ln\left(\frac1{|x-y|}\right)\right)F(u(y))\dd y\\
		&=\intd\ln(1+|y|)F(u(y))\dd y<\|F(u)\|_{\ast,1}<+\infty\,,
	\end{aligned}
\end{equation*}
that is constant. This implies that $\Phi_u$ is a solution of \eqref{SP:system} in the sense of Definition \ref{sol_SP}, by applying \cite[Lemma 2.4]{H}, for which all such solutions of $-\Delta\Phi=\ff$ in $\R^2$ are of the form $\Phi=\tw_u+p$ with $p$ polynomial of degree at most $1$.

\vskip0.4truecm
\paragraph{\textbf{Acknowledgements}:} The main part of this work was carried out while F.B. was affiliated with Università degli Studi di Milano and G.R. with Università degli Studi dell'Insubria, whose support is gratefully acknowledged. The authors are members of \textit{Gruppo Nazionale per l'Analisi Matematica, la Probabilità e le loro Applicazioni} (GNAMPA) of the \textit{Istituto Nazionale di Alta Matematica} (INdAM). Federico Bernini is partially supported by INdAM-GNAMPA Project 2024 titled \textit{Problemi spettrali e di ottimizzazione di forma: aspetti qualitativi e stime quantitative}. Giulio Romani and Cristina Tarsi are partially supported by INdAM-GNAMPA Project 2023 titled \textit{Interplay between parabolic and elliptic PDEs} and INdAM - GNAMPA Project 2024 titled \textit{New perspectives on Choquard equation through PDEs with local sources}.

\bibliographystyle{abbrv} 
\bibliography{Bibliography.bib}

\end{document}